\documentclass[a4paper,10pt]{article}
\pdfoutput=1
\usepackage[scaled]{helvet}              % use font Helvetica
\usepackage[T1]{fontenc}                 % for support of Helvetica
\usepackage{geometry}
\usepackage{setspace}                 % for line spacing
\usepackage{exscale}                  % pack­age for scaling of the math extension font cmex
\usepackage{a4wide}                   % for a4 format
\usepackage{framed}                   % package for frames
\usepackage{subscript}                % provides the commands \textsuperscript and \textsubscript
\usepackage{authblk}                  % for listing of authors
\usepackage{appendix}                 % special numbering for the appendix. Use \begin{appendix} ... \end{appendix}.
\usepackage{colordvi}                 % provides more colors
\usepackage{enumitem}                 % provides list environment \begin{enumerate} \end{enumerate} with various additional features
\usepackage[american]{babel}
\usepackage{fancyhdr}                   %package for designing head and foot
\usepackage[colorlinks=true,linkcolor=black,citecolor=black,filecolor=black,urlcolor=black,hypertexnames=false]{hyperref}
\usepackage{aliascnt}      % for defining clever refs
\usepackage{cleveref}      % clever refences. include AFTER hyperref and BEFORE amsmath and ntheorem
\usepackage[babel,german=guillemets]{csquotes}
\usepackage[backend=bibtex,style=numeric-comp,useprefix=true,hyperref=true,firstinits=true]{biblatex}
\bibliography{references.bib}
\usepackage[fleqn]{amsmath} % math symbols
\usepackage{mathtools}      % math symbols
\usepackage{amssymb}        % math symbols
\usepackage{latexsym}       % provides 11 special math symbols
\usepackage{dsfont}         %provides symobol for unity matrix \mathds{1}
\numberwithin{equation}{section}      %equation numbering subordinate to chapters\usepackage{enumerate}
\usepackage{amsthm}         % math theorem environments

\usepackage[framemethod=PSTricks]{mdframed}

\theoremstyle{plain}    
   \crefformat{thm}{Theorem~#2#1#3}                %define citation format for counter "thm"
   \newtheorem{thm}{Theorem}[section]              %define new environment "thm" that automatically creates a counter "thm". Choose numbering subordinate to chapters    
   \crefname{thm}{Theorem}{Theorems}               %define single/plural citation names for environment "thm"
   \crefformat{count_prop}{Proposition~#2#1#3}        %define citation format for counter "count_prop"
   \newaliascnt{count_prop}{thm}                      %define counter "count_prop" as a alias for the counter "thm". This enables consecutive numbering.
   \aliascntresetthe{count_prop}                      %close definition
   \crefname{proposition}{Proposition}{Propositions}  %define single/plural citation names for environment "proposition"
   \crefformat{count_lem}{Lemma~#2#1#3}       
   \newaliascnt{count_lem}{thm}               
   \newtheorem{lemma}[count_lem]{Lemma}        
   \aliascntresetthe{count_lem}                
   \crefname{lemma}{Lemma}{Lemmas}
   \crefformat{count_defi}{Definition~#2#1#3}       
   \newaliascnt{count_defi}{thm}
   \newtheorem{definition}[count_defi]{Definition}
   \aliascntresetthe{count_defi}                
   \crefname{definition}{Definition}{Definitions}
   \crefformat{count_cor}{Corollary~#2#1#3}       
   \newaliascnt{count_cor}{thm}
    
   \aliascntresetthe{count_cor}                
   \crefname{corollary}{Corollary}{Corollaries}

\theoremstyle{definition}  
   \crefformat{count_rem}{Remark~#2#1#3}       
   \newaliascnt{count_rem}{thm}
   \newtheorem{remark}[count_rem]{Remark}
   \aliascntresetthe{count_rem}                
   \crefname{remark}{Remark}{Remarks}
   \crefformat{count_exmpl}{Example~#2#1#3}       
   \newaliascnt{count_exmpl}{thm}
   
   \aliascntresetthe{count_exmpl}                
   \crefname{example}{Example}{Examples}

   \crefformat{fthm}{Theorem~#2#1#3}               %define citation format for counter "fthm"
   \mdtheorem{fthm}{Theorem}[section]              %define new environment "fthm" that automatically creates a counter "fthm". Choose numbering subordinate to chapters    
   \crefname{fthm}{Theorem}{Theorems}              %define single/plural citation names for environment "thm"
   \crefformat{count_fprop}{Proposition~#2#1#3}        %define citation format for counter "count_fprop"
   \newaliascnt{count_fprop}{fthm}                     %define counter "count_fprop" as a alias for the counter "fthm". This enables consecutive numbering.
   \mdtheorem{fproposition}[count_fprop]{Proposition}  %define new environment "fproposition" that uses the counter "count_fprop" with the above specified ciation format 
   \aliascntresetthe{count_fprop}                      %close definition
   \crefname{fproposition}{Proposition}{Propositions}  %define single/plural citation names for environment "fproposition"
   \crefformat{count_flem}{Lemma~#2#1#3} 
   \newaliascnt{count_flem}{fthm}
   \mdtheorem{flemma}[count_flem]{Lemma}                
   \aliascntresetthe{count_flem}
   \crefname{flemma}{Lemma}{Lemmas}
   \crefformat{count_fdefi}{Definition~#2#1#3}
   \newaliascnt{count_fdefi}{fthm}
   \mdtheorem{fdefinition}[count_fdefi]{Definition}      
   \aliascntresetthe{count_fdefi}
   \crefname{fdefinition}{Definition}{Definitions}
   \crefformat{count_fcor}{Corollary~#2#1#3}
   \newaliascnt{count_fcor}{fthm}
   \mdtheorem{fcorollary}[count_fcor]{Corollary}       
   \aliascntresetthe{count_fcor}
   \crefname{fcorollary}{Corollary}{Corollary}
   \crefformat{count_frem}{Remark~#2#1#3}
   \newaliascnt{count_frem}{fthm}
   \mdtheorem{fremark}[count_frem]{Remark}
   \aliascntresetthe{count_frem}
   \crefname{fremark}{Remark}{Remarks}
   \crefformat{count_fexmpl}{Example~#2#1#3}       
   \newaliascnt{count_fexmpl}{thm}
   \mdtheorem{fexample}[count_fexmpl]{Example}
   \aliascntresetthe{count_fexmpl}                
   \crefname{fexample}{Example}{Examples}

\usepackage{xargs}   %usepackage for macros with several default parameters
\newcommand{\GagNirenberg}{Gagliardo-Nirenberg's inequality}
\newcommand{\Grownwall}{Gronwall's inequality}      % catch misspelling
\newcommand{\Holder}{H\"{o}lder's inequality}

\newcommand{\Young}{Young's inequality}
\newcommand*{\pnp}{Poisson--Nernst--Planck system}
\newcommand*{\dpnp}{Darcy--Poisson--Nernst--Planck system}

\newcommand{\Div}{\ensuremath{ \textnormal{div} }}                          % \Div
\newcommand{\Sign}{\ensuremath{ \textnormal{sign} }}                         % \Sign
\newcommand{\weak}{\ensuremath{ \textnormal{weak} }}                        % \weak
\newcommand{\weakstar}{ \textnormal{weak}^\ast }                            % \weakstar
\newcommand{\vphi}{\ensuremath{ \varphi }}

\newcommand{\setN}{\ensuremath{ \mathbb{N} }}                   
\newcommand{\setZ}{\ensuremath{ \mathbb{Z} }}                   
\newcommand{\setR}{\ensuremath{ \mathbb{R} }}                   
                                                               
\newcommand*{\OmegaT}[1][T]{\ensuremath{ \Omega_{#1} }}      % parabolic cylinder: \OmegaT[x] 
\newcommand*{\GammaT}[1][T]{\ensuremath{ \Gamma_{#1} }}      % parabolic boundary: \GammaT[x]
\newcommand*{\brac}[1]{\ensuremath{ \left( {#1} \right) }}                       % bracket: \brac{xx}                     
\newcommand*{\sqbrac}[1]{\ensuremath{ \left[ {#1} \right] }}                     % squared bracket: \sqbrac{xx}                     
\newcommand*{\cbrac}[1]{\ensuremath{ \left\{ {#1} \right\} }}                    % curly bracket \cbrac{xx}                       
\newcommand*{\Hence}{\ensuremath{ ~~\Longrightarrow~~ }}                         % long implication arrow with more spacing
\newcommand*{\Equivalent}{\ensuremath{ ~~\Longleftrightarrow~~ }}                % long equivalence arrow with more spacing   
\renewcommand*{\vec}[1]{ \ensuremath{ \boldsymbol{#1} }}  % Vector: \vec{x}
\newcommand*{\abs}[1]{\left| {#1} \right| }               % Vectornorm: \abs{x} 
\newcommand*{\vecA}{\vec{A}}

\newcommand*{\vecc}{\vec{c}}

\newcommand*{\vecE}{\vec{E}}

\newcommand*{\vecf}{\vec{f}}

\newcommand*{\vecnu}{\vec{\nu}}
\newcommand*{\vecnull}{\vec{0}}

\newcommand*{\vecsigma}{\vec{\sigma}}

\newcommand*{\vecu}{\vec{u}}
\newcommand*{\vecv}{\vec{v}}

\newcommand*{\scp}[3][~]{\ensuremath{ \left( {#2}{#1},{#1}{#3} \right) }}    % scalar product:  \scp[spacing]{x}{y} 
\newcommand*{\dualp}[3][~]{\ensuremath{ \left\langle {#2}{#1},{#1}{#3} \right\rangle }} % dual pairing:   \dualp[spacing]{x}{y}        
\newcommand*{\grad}{\vec{\nabla}\!}                       % gradient: \grad 
\newcommand*{\dert}[1][t]{\ensuremath{ \partial_{{#1}} }} % temporal derivative: \dert[opt]   
\newcommand*{\derr}[1][t]{\ensuremath{ \frac{d}{d{#1}} }} % dervative on real line: \derr[opt]
\newcommand{\Int}[4]{\ensuremath{ \int\limits_{#2}^{{#3}} {#4} d{#1} }}    % Integral:   \Int{integrand variable}{subscript}{superscript}{function}  
\newcommand*{\Intdx}[2][\Omega]{ \Int{x}{{#1}}{}{{#2}} }                   % Volume integral with respect to dx and defaluft domain: \Intdx[domain]{function}
\newcommandx*{\Intdt}[3][1=0,2=T,usedefault]{ \Int{t}{{#1}}{{#2}}{{#3}} }  % Time integral with respect to dt and default interval:  \Intdt[lower border][upper border]{function}
\newcommandx*{\fspace}[4][1=-1pt,usedefault]{ \ensuremath{  {{#2}} \hspace{#1} \left( {#3}{#4} \right) }} % Function space: \fspace[hspace]{name}{domain}{range} 
\newcommandx*{\norm}[2]{ \ensuremath{ \left\| {#1} \right\|_{{#2}}   }}                                   % Norms for functions: \norm{object}{space} 
\newcommandx*{\Ck}[4][1=\Omega,2= ,usedefault]{ \fspace{ C^{{#3}}_{{#2}} }{ {#1} }{ {#4} } } %\Ck[domain][subscript]{smoothness k}{range}  
\newcommandx*{\Lp}[4][1=\Omega,2= ,usedefault]{ \fspace{ L^{{#3}}_{{#2}} }{ {#1} }{ {#4} } } %\Lpk[domain][subscript]{exponent p}{range} 
\newcommandx*{\Wkp}[5][1=\Omega,2= ,usedefault]{ \fspace{ W^{{#3},{#4}}_{{#2}} }{ {#1} }{ {#5} } } % \Wkp[domain][subscript]{smoothness k}{exponent p}{range}                
\newcommandx*{\Hk}[4][1=\Omega,2= ,usedefault]{ \fspace{ H^{{#3}}_{{#2}} }{ {#1} }{ {#4} } }       % \Hk[domain][subscript]{smoothness k}{range}                      
\newcommandx*{\Hkdiv}[3][1=\Omega,2= ,usedefault]{ \fspace{ H^{{#3}}_{{#2}} }{\Div}{;{#1}} }       % \Hkdiv[domain][subscript]{smoothness k}

\usepackage{xargs}   %usepackage for macros with several default parameters 
\newcommand{\spaceTC}{ \Hk{1}{} }                                                         % Concentrations: Test function space
\newcommand{\spaceSC}{ \fspace{L^2}{I}{; \Hk{1}{} } \cap \fspace{H^1}{I}{; \Hk{1}{}^* } } % Concentrations: Solution space
\newcommand{\spaceTF}{ \Hkdiv[][0]{} }                                                    % Velocity: Test function space
\newcommand{\spaceSF}{ \fspace{L^\infty}{I}{; \Hkdiv[][f]{} } }                           % Velocity: Solution space
\newcommand{\spaceTP}{ \Lp{2}{} }                                                         % Pressure: Test function space
\newcommand{\spaceSP}{ \fspace{L^\infty}{I}{; \Lp{2}{}/\setR } }                          % Pressure: Solution space
\newcommand{\spaceSEF}{ \fspace{L^\infty}{I}{; \Hkdiv[][\sigma]{} }  }                    % Electrostatic Potential: Solution space
\newcommand*{\tx}{(t,x)}    % continuous variables (t,x)
\newcommand*{\timeEnd}{T_0} % end time T_0

\newcommand*{\test}{\ensuremath{ \vphi }} % Test function: scalar
\newcommand*{\Test}{\ensuremath{ \vecv }} % Test function: vector

\newcommand*{\boltz}{\ensuremath{ k_b }}                                                 % Parameters: Boltzmann Constant
\newcommand*{\temp}{\ensuremath{ T }}                                                    % Parameters: Temperature
\newcommand*{\Dl}[1][]{\ensuremath{ \mathcal{D}_{{#1}} }}                                % Parameters: Diffusion coefficient
\newcommand*{\zl}[1][l]{\ensuremath{ z_{{#1}} }}                                         % Parameters: Valency
\newcommand*{\mob}[1][]{\ensuremath{ \omega_{{#1}} }}                                    % Parameters: Mobility
\newcommand*{\wl}[1][l]{\ensuremath{ e\zl[{#1}]\mob }}                                   % Parameters: Electrophoretic mobility
\newcommand*{\permitELscalar}{\ensuremath{ \epsilon }}                                   % Parameters: Electric permittivity constant
\newcommand*{\permitEL}{\ensuremath{ \mathcal{E} }}                                      % Parameters: Electric permittivity tensor
\newcommand*{\permeabH}{\ensuremath{ \mathcal{K} }}                                      % Parameters: Hydraulic permeability
\newcommand*{\coeffEL}[1][l]{\ensuremath{e\zl[{#1}](\permitELscalar\boltz\temp)^{-1} }}  % Parameters: Electric coefficient in Transport
\newcommand*{\coeffForceEL}{\ensuremath{\mu^{-1}\permitEL^{-1}}}

\newcommand{\cl}[1][l]{\ensuremath{ c_{{#1}} }}                  % Concentrations: l-th concentration
\newcommand{\Cl}{\ensuremath{ \vecc }}                           % Concentrations: vector
\newcommand{\clstart}[1][l]{\ensuremath{ c_{0,{#1}} }}           % Concentrations: Initial value for l-th concentration
\newcommand{\gl}[1][l]{ g_{{#1}} }                               % Concentrations: boundary flux of l-th concentration
\newcommand{\Rl}[1][l]{\ensuremath{ R_{{#1}} }}                                 % Reactions: reaction rate in equation for l-th concentration   
\newcommand{\fieldF}{\ensuremath{ \vecu }} % Velocity: Flow field
\newcommand{\press}{p}                     % Pressure: Pressure field

\newcommand{\potEL}{\ensuremath{ \Phi }}                                    % Electrostatic:  Electric Potential
\newcommand{\fieldEL}{\ensuremath{ \vecE }}                                 % Electrostatic:  Electric field
\newcommand*{\chargeEL}[1][]{ \ensuremath{ \rho_{f{#1}} } }                 % Electrostatic:  Free charge density -- short
\newcommand*{\chargeELlong}[1][l]{ \sum_l\zl\cl }                           % Electrostatic:  Free charge density -- long
\newcommand*{\chargeELlongTwo}[1][l]{ \zl[1]\cl[1]-\abs{\zl[2]}\cl[2] }     % Electrostatic:  Free charge density 2 species -- long
\newcommand*{\forceEL}[1][]{\ensuremath{ \chargeEL[{#1}] \fieldEL_{{#1}} }} % Electrostatic:  Electroosmotic force -- short
\newcommand{\clm}[1][l]{\ensuremath{ c_{{#1}}^m }}        % Upper Truncated Concentrations: l-th concentration
\newcommand{\fixOp}{\ensuremath{\mathcal{F}}}
\newcommand{\fixSet}{\ensuremath{K}}

\newcommand{\clold}[1][l]{\ensuremath{\bar{c}_{{#1}} }}                           % Old Concentrations
\newcommand{\Clold}{\ensuremath{ \bar{\vecc} }}                                   % Old Concentration vector
\newcommand*{\chargeELold}[1][]{ \ensuremath{ \bar{\rho}_{f{#1}} } }              % Old free charge density -- short
\newcommand*{\chargeELlongold}[1][l]{ \sum_l\zl\clold }                           % Old free charge density -- long
\newcommand*{\chargeELlongTwoold}[1][l]{ \zl[1]\clold[1]-\abs{\zl[2]}\clold[2] }  % Old free charge density 2 species -- long
\newcommand*{\forceELold}[1][]{\ensuremath{ \chargeELold[{#1}] \fieldEL_{{#1}} }} % Old electroosmotic force -- short
\crefname{figure}{Figure}{Figures}
\crefname{algorithm}{Algorithm}{Algorithms}
\crefname{subsection}{Section}{Section}
\crefname{section}{Section}{Section}
\crefname{chapter}{Chapter}{Chapter}
\crefname{part}{Part}{Part}
\crefname{appendix}{Appendix}{Appendices}
\newcounter{countAssumption} % counter for consecutive numbering of assumption lists
\newcounter{countNotation}   % counter for consecutive numbering of notation lists
\newcommand{\labelInt}{I.a}
\begin{document}
%
%%%%%%%%%%%%%%%%%%%%%%
% 1.5 line spacing
%%%%%%%%%%%%%%%%%%%%%%
% \onehalfspacing
%
%%%%%%%%%%%%%%%%%%%%%%%%%%%%%
%  TITLE                    %
%%%%%%%%%%%%%%%%%%%%%%%%%%%%%
\title{Global existence of weak solutions of a model for electrolyte solutions -- Part 1: Two-component case}
\author[1]{Matthias Herz}
\author[1]{Peter Knabner}
\renewcommand\Affilfont{\itshape\small}
\affil[1]{Department of Mathematics, University of Erlangen-N\"urnberg, Cauerstr. 11, D-91058 Erlangen, Germany}
\date{\today}
\maketitle
%
%
%%%%%%%%%%%%%%%%%%%%%%%%%%%%%
%  ABSTRACT                 %
%%%%%%%%%%%%%%%%%%%%%%%%%%%%%
\begin{abstract}
This paper analytically investigates the \dpnp. This system is a mathematical model for electrolyte solutions. In this paper, we consider electrolyte solutions,
which consist of a neutral fluid and two suspended oppositely charged chemical species with arbitrary valencies~$\zl[1]>0>\zl[2]$. We prove global existence and uniqueness of
weak solutions in two space dimensions and three space dimensions.
\par
So far, most of the existence results have been proven for symmetric electrolyte solutions. These solutions consist of a neutral fluid and two suspended charged chemical species
with symmetric valencies~$\pm\zl[]$. As many electrolyte solutions in biological applications and hydrodynamical applications are not symmetric, the presented extension of the
previous existence results is an important step.
\\[2.0mm]
\end{abstract}
%
%
%%%%%%%%%%%%%%%%%%%%%%%%%%%%%
%%%%%%%%%%%%%%%%%%%%%%%%%%%%%
%  SECTION introduction     %
%%%%%%%%%%%%%%%%%%%%%%%%%%%%%
%%%%%%%%%%%%%%%%%%%%%%%%%%%%%
\section{Introduction}\label{sec:Introduction}
%% Motivation
Many complicated phenomena in hydrodynamics and biology can be modeled in the context of electrolyte solutions.
The reason for this is that models for electrolyte solutions must simultaneously capture the following three ubiquitous processes:
(i) the transport of the charged particles, (ii) the hydrodynamic fluid flow, (iii) the electrostatics.
Moreover, these processes simultaneously occur in electrolyte solutions. Firstly, the electrostatic field is generated by the movement of the charged particles,
and conversely, the movement of the charged particles is influenced by the electrostatic field. Secondly, the fluid flow changes the flux of the charged particles
and conversely, the moving charged particles lead to a force term, which generates an electroosmotic fluid flow.
\medskip
\par
%% General overview
The classical models for electrolyte solutions, that capture the fully coupled nature of these processes are the so-called \pnp{s} (for a fluid at rest) and the \dpnp{s} (for laminar flow in porous media).
In particular, \pnp{s} are also known as drift--diffusion systems, van Rosenbroeck systems, or semiconductor device equations.
For a detailed derivation of these systems, we refer to \cite{DreyerGuhlkeMuller, Elimelech-book, Juengel-book, Masliyah-book, Russel-book, Samohyl, Schuss, Mielke10, Probstein-book}.
Among many others, these models have been investigated analytically in
\cite{BT2011, Burger12, Castellanos-book, GasserJungel97, Glitzky2004, Gajewski, HyongEtAll, Frehse-book, Markovich-book, Mielke10, Roubicek2005-1, Schuss, Wolfram, Schmuck11, herz_ex_dnpp}.
\medskip
\par
%% History
So far, most of the analytical investigations have been carried out for electrolyte solutions, which consist of a electrically neutral solvent (at rest) and two oppositely charged chemical species with
symmetric valencies~$\pm z$. One reason for this is that especially the symmetric valencies~$\pm1$ naturally occur in the context of semiconductor devices and most of the analytical investigations
are related to semiconductor devices. Previous existence results, which consider charged solutes with arbitrary valencies, were proven amongst others in
\cite{BT2011, Burger12, Glitzky2004, Roubicek2005-1, Roubicek2006}. These papers considered electrolyte solutions with multiple suspended charged solutes. We investigate this multicomponent case
in the second part of this work, whereas in this paper, we focus on the two-component case. This means, we consider electrolyte solutions, that consist of a neutral solvent and two oppositely charged solutes.
In this situation, the results of this paper go beyond the above mentioned papers. More precisely, the authors of \cite{BT2011} proved local in time existence. The results of \cite{Roubicek2005-1, Roubicek2006}
were proven under the additional assumption of a volume-additivity constraint and by including an additional reaction force term in the transport equations. These additional assumptions allow to
bypass in \cite{Roubicek2005-1, Roubicek2006} the main difficulties, which we briefly sketch below. Finally, the paper~\cite{Burger12} dealt with a stationary system and existence in two dimensions was
established in \cite{Glitzky2004}.
\medskip
\par
%% Explanation of the new Problems
In the proof of the crucial a~priori estimates occur additional difficulties, if we allow for two oppositely charged solutes with arbitrary valencies~$\zl[1]>0>\zl[2]$. More precisely, the main difficulty is
to obtain a~priori estimates for the solutes~$\cl$, which are independent of the electric field. Such a~priori estimates are easily obtained in case of symmetric valencies~$\pm z$. To briefly sketch this,
we consider two charged solutes~$\cl[1]$ (positively charged) and $\cl[2]$ (negatively charged) with symmetric valencies~$\pm z$. In the proof of a~priori estimates for $\cl$, we test the
equations for $\cl$ with the standard test functions~$\test=\cl$, and we remember that the electric field~$\fieldEL$ satisfies according to Gauss's law~$\grad\cdot\vecE=z(\cl[1]-\cl[2])$. Thereby, we obtain for
the sum of the \enquote{electric~drift integrals}, which describe the electrophoretic motion
\begin{align*}
 &-2\sum_l \scp{\pm z\cl\vecE}{\grad\cl~}_{\Lp{2}{}} =z\scp{\grad\cdot\vecE}{(\cl[1])^2-(\cl[2])^2}_{\Lp{2}{}} \\
 &= z^2 \scp{\cl[1] - \cl[2]}{(\cl[1])^2 -(\cl[2])^2}_{\Lp{2}{}}  \geq 0, \quad \text{since } \sqbrac{a - b}\sqbrac{a^2 - b^2} \geq 0 ~\text{ for } a,b\geq0~.
\end{align*}
Due to this pointwise sign condition, we can omit the sum of the \enquote{electric~drift integrals} and the a~priori estimates for $\cl$ are naturally independent of the electric field~$\fieldEL$.
In case of arbitrary valencies, such a pointwise sign condition does not hold true, if we use the standard test functions~$\test=\cl$. However, we propose to carry over this pointwise sign condition
by using weighted test functions~$\test=\abs{\zl}\cl$ instead.
\medskip
\par
% Contribution
The contribution of this paper is to prove global existence, uniqueness, and boundedness of weak solutions for electrolyte solutions, which consist of a neutral solvent and
two charged solutes with arbitrary valencies. Moreover, we do not impose any further restrictions such as the often used electroneutrality constraint, cf. \cite{Allaire10},
or the volume additivity constraint, cf. \cite{Roubicek2005-1}. This result is a first step towards the treatment of multicomponent electrolyte solutions, which contain $L\in\setN$ solutes.
The presented proof is based on the weighted test functions~$\test=\abs{\zl}\cl$. More precisely, we particularly use these weighted test functions in the proof of \cref{lemma:energy},
which is the basis for the following a~priori estimates in \cref{lemma:bounded} and \cref{thm:aprioriBounds}.
\medskip
\par
%% Structure of the paper
The rest of this paper is organized as follows: In \cref{sec:Model}, we present the \dpnp\ and in \cref{sec:Uniqueness}, we prove that solutions are unique.
Then, we introduce the fixed point operator in \cref{sec:fpi-operator}. Finally, we show the crucial a~priori estimates in \cref{subsec:aprioriBounds}, and in \cref{subsec:fpi-point},
we show the global existence.
%
%
%%%%%%%%%%%%%%%%%%%%%%%%%%%%%
%%%%%%%%%%%%%%%%%%%%%%%%%%%%%
%  SECTION model equations
%%%%%%%%%%%%%%%%%%%%%%%%%%%%%
%%%%%%%%%%%%%%%%%%%%%%%%%%%%%
\section{Model Equations}\label{sec:Model}
Subsequently, we present the \dpnp. This system is a field-scale model%
\footnote{For a detailed introduction to the modeling of porous media and the notion of field-scales and pore-scales, we refer to  \cite{bear-book}.} %
for electrolyte solutions in porous media. A rigorous derivation of field-scale \dpnp{s} from pore-scale systems was carried out, e.g., in \cite{Allaire10, RayMunteanKnabner}.
Note that commonly on field-scales volume effects dominate and surface effects such as the electrostatic double-layer effects are negligible.
However, a characteristic feature of porous media is are dominating surface effects even on field scales. This justifies to consider field-scale \dpnp{s}.
\par
We now introduce some notation, in order to present the model equations.
%=============================
% LIST: NOTATION
%=============================
\begin{enumerate}[label=({N}\arabic*), ref=({N}\arabic*), itemsep=1mm, start=\value{countNotation}+1]
 \item \textbf{Geometry: } For $n\in\{2,3\}$, let $\Omega\in\setR^n$ be a $n$-dimensional bounded domain with boundary~$\Gamma$ and
       corresponding exterior normal field~$\vecnu$.
       Next, let $I:=(0,\timeEnd)$ be a time interval and we introduce by $\OmegaT:= I\times\Omega$
       a time space cylinder with lateral boundary~$\GammaT:=I\times\partial\Omega$.
       Furthermore, we suppose $\Omega$ to be a porous medium with constant porosity~$\theta$.%
       \label{Not:Geom}
 \item \textbf{Variables: } We assume that $\Omega$ is fully saturated with a fluid, in which two charged chemical species are suspended.
       We denote the velocity field of the electrolyte solution by $\fieldF$, its pressure by $p$, the electric field and the electrostatic potential
       in the electrolyte solution by $\vecE$ and $\potEL$. Next, we denote the number densities of the respective chemical species by~$\cl$,
       $l\in\{1,2\}$. Furthermore, we define the concentration vector by $\vecc:=(\cl[1],\cl[2])$.%
       \label{Not:Variables}
 \item \textbf{Electrics: } The chemical species~$\cl$ carry a charge~$e\zl$. Here, $e$ is the elementary charge
       and $\zl\in\setZ$ the respective valency ($\zl\neq0$). W.l.o.g., we assume $\zl[1]>0>\zl[2]$. The chemical species~$\cl$ possess electric mobilities~$\wl$, where $\mob$ is
       the so-called mobility tensor. It is~$\Dl=\boltz\temp\mob$ according to Einstein-Smoluchowski relation, see \cite[Chapter~6]{Masliyah-book}.
       Here, $\boltz$ is the Boltzmann constant, $\temp$ the temperature. Hence, we have the identity $\wl=e\zl(\boltz\temp)^{-1}\Dl$.
       We denote by $\chargeEL$ the free charge density and by $\rho_b$ a background charge density, e.g.,
       coming from not resolved pore-scale inclusions inside $\Omega$.%
       \label{Not:Electrics}
 \item \textbf{Coefficients: } We denote by $\Dl$ the diffusion-dispersion tensor, which is identical for all chemical species~$\cl$.
       Although the molecular diffusion might be different for each $\cl$, the dispersion coming from a tortuous geometry
       is by far dominating on the considered field scales. Since the geometry looks the same for all chemical species,
       we obtain a coinciding diffusion-dispersion tensor~$\Dl$. Next, we denote by $\permeabH$ the constant
       permeability tensor of the medium, by $\mu$ the dynamic viscosity of the fluid, and by $\permitEL:=\permitELscalar\Dl$
       the constant electric permittivity tensor of the medium. For a rigorous derivation of the last relation, see \cite{RayMunteanKnabner}.
       We note, that we have the identity $\permitEL =\permitELscalar\Dl=\permitELscalar\boltz\temp\mob$.%
       \label{Not:coeff}
\setcounter{countNotation}{\value{enumi}}
\end{enumerate}
We suppose the boundary~$\Gamma$ of the domain~$\Omega$ is charged, e.g., from surfactants. As the solutes~$\cl$ carry charges,
they interact with $\Gamma$ in a small boundary layer, the so-called electrostatic double-layer. This leads to a spatially inhomogeneous charge distributions,
which gives rise to an electric field~$\fieldEL$. Simultaneously, the electric field~$\fieldEL$ generates an electric body force in the surrounding fluid.
Thereby, an electroosmotic flow develops, which in turn interacts with the chemical species. This leads to an interplay between the
electrophoretic movement of the charged particles, the electroosmotic flow of the fluid, and a varying electric field.
\par
\dpnp{s} capture these coupled processes based on the following three conservation laws:\\[2.0mm]
\begin{subequations}
%
%=============================
% MODEL: Gauss's Law
%=============================
\textbf{Law 1 -- Gauss's law:} The surface charges and the charged solutes~$\cl$
give rise to an electric field~$\fieldEL$. For the electric field~$\fieldEL$, we solve Gauss's law.
Additionally, we assume that the electric field is generated by an electrostatic potential~$\potEL$.
Thus, we have $\vecE:=-\grad\potEL$. The boundary data are denoted by $\sigma$ and the initial data
are obtained by substituting the initial data~$\clstart$ of the charged solutes~$\cl$ on the right hand side
of Gauss's law and solving for the electric field~$\fieldEL$. Furthermore, we assume that inside the electrolyte
solution, we have a background charge density~$\rho_b$, coming, e.g., from not resolved pore-scale inclusions inside~$\Omega$.
Mathematically, Gauss's law writes for the redefined electric field $\vecE:=\permitEL\fieldEL$ as
\vspace{-2mm}
\begin{align}
              \permitEL^{-1}\fieldEL   &=  -\grad\potEL                &&  \text{in } \OmegaT,   \label{eq:Model-1a} \\
           \grad\cdot\fieldEL          &=   \chargeEL + \rho_b         &&  \text{in } \OmegaT,   \label{eq:Model-1b} \\
                           \chargeEL   &=  \theta(\chargeELlongTwo)    &&  \text{in } \OmegaT,   \label{eq:Model-1c} \\
             \fieldEL\cdot\vecnu       &=  \sigma                      &&  \text{on } \GammaT,   \label{eq:Model-1d}
\end{align}
%
%=============================
% MODEL: Darcy's Law
%=============================
\textbf{Law 2 -- Darcy's law:} The velocity field~$\fieldF$ is subject to conservation of mass and momentum.
On field-scales, this is sufficiently well-captured by Darcy's law, which connects the velocity field~$\fieldF$ and the
pressure gradient~$\grad\press$. As we include electroosmotic flows, an electric body force term enters the equations.
The boundary data are denoted by $f$ and the initial data are obtained by inserting $\clstart$ and $\fieldEL(0)$ on the right hand side of Darcy's law.
Mathematically, Darcy's law reads (with the redefined electric field $\vecE:=\permitEL\fieldEL$) as
\vspace{-2mm}
\begin{align}
    \permeabH^{-1} \fieldF    &= \mu^{-1}\brac{-\grad\press + \permitEL^{-1}\forceEL}  &&  \text{in } \OmegaT,  \label{eq:Model-1e}  \\
	  \grad\cdot\fieldF   &=  0                                                    &&  \text{in } \OmegaT,  \label{eq:Model-1f}  \\
	  \fieldF\cdot\vecnu  &=  f                                                    &&  \text{on } \GammaT.  \label{eq:Model-1g}
\end{align}
%
%=============================
% MODEL: Nernst-Planck PDEs
%=============================
\textbf{Law 3 -- Nernst--Planck equations:} The evolution of the chemical species~$\cl$ is subject to mass continuity.
Here, the mass flux arises due to diffusion, convection, and an electric~drift. Such mass fluxes are called Nernst--Planck fluxes.
We assume the equations for $\cl$ are coupled through reaction rates~$\Rl$.
The initial data are denoted by $\clstart$ and the flux boundary data by~$\gl$.
Mathematically, Nernst--Planck equations are given (with \ref{Not:coeff} and the redefined electric field $\vecE:=\permitEL\fieldEL$) by
\vspace{-2mm}
\begin{flalign}
  \theta\dert\cl + \grad\cdot\brac{\Dl\grad\cl + {\cl}[\fieldF + \coeffEL\fieldEL] } &= \theta\Rl(\Cl)   &&  \text{in } \OmegaT,  & \label{eq:Model-1h} \\[0.2cm]
                  \brac{\Dl\grad\cl + {\cl}[\fieldF + \coeffEL\fieldEL] }\cdot\vecnu &= \gl              &&  \text{on } \GammaT,  & \label{eq:Model-1i} \\[0.2cm]
                                                                              \cl(0) &= \clstart         &&  \text{on } \Omega.   & \label{eq:Model-1j}
\end{flalign}
\end{subequations}
\begin{remark}
\textnormal{
Equations~\eqref{eq:Model-1a}--\eqref{eq:Model-1d} are Poisson's equation for $\potEL$ in mixed formulation.
For this reason \enquote{Poisson} is contained in the name \dpnp.
For analytical investigations, it is of advantage to deal with Poisson's equation directly, as the comprehensive regularity results for Poisson's equation hold true, cf. \cite{GilbargTrudinger-book}.
However, the mixed formulation is of advantage, as we can easily introduce in \eqref{eq:Model-1a} the general electric fields~$\fieldEL= -\grad\potEL -\dert\vecA$,
which is the expression of an electric field according to Maxwell's equations in terms of the electromagnetic potentials, cf. \cite{LifshitzLandau-book2}.
Furthermore, the mixed formulation is of advantage as starting point for numerical approximations, since this leads to a direct approximation of the electric field~$\fieldEL$, cf. \cite{FrankRay2011_pnp}.%
}\hfill$\square$
\end{remark}
\begin{remark}
\textnormal{
The boundary flux in equation~\eqref{eq:Model-1i} can be equivalently expressed with equations \eqref{eq:Model-1d} and \eqref{eq:Model-1g} by
\begin{align*}
 \brac{\Dl\grad\cl}\cdot\vecnu &= \gl - {\cl}f - (\permitELscalar\boltz\temp)^{-1}e\zl{\cl}\sigma  \qquad \text{on }~ \GammaT~.
\end{align*}
Thus, the boundary flux condition is equivalent to a Robin boundary condition for the diffusion part.
}\hfill$\square$
\end{remark}
\begin{remark}
\textnormal{
Equations~\eqref{eq:Model-1a}--\eqref{eq:Model-1j} contain the nonlinear coupling terms $\forceEL$ in Darcy's law,
and $\cl\fieldF$, $\cl\fieldEL$ in Nernst--Planck equations. These nonlinearities arise only after combing the three subsystems to a \dpnp.
This reflects the fact, that the coupling of initially isolated subprocesses leads to new nonlinearities in the resulting system.%
}\hfill$\square$
\end{remark}
%
%%%%%%%%%%%%%%%%%%%%%%%%%%%%%%%%%%%%%%%%%%%%%%%%%%%%
% SUBSECTION notation Assumption, Weak Formulation %
%%%%%%%%%%%%%%%%%%%%%%%%%%%%%%%%%%%%%%%%%%%%%%%%%%%%
\subsection{Notation, Assumptions and Weak Formulation}
We now introduce the required notation for the analytical investigations.
%
%=============================
% LIST: Notation
%=============================
\begin{enumerate}[label=({N}\arabic*), ref=({N}\arabic*), itemsep=0.0mm, start=\value{countNotation}+1]
 \item \textbf{Spaces: } For $k>0$, $p\in [1,\infty]$, we denote the Lebesgue spaces for scalar-valued and vector-valued functions by $\Lp{p}{}$
       and the respective Sobolev spaces by $\Wkp{k}{p}{}$, cf.~\cite{Adams2-book}.
       Furthermore, we set $\Hk{k}{}:=\Wkp{k}{2}{}$ and we refer for the definition of the Bochner spaces~$\fspace{L^p}{I}{;V}$, $\fspace{H^k}{I}{;V}$
       over a Banach space~$V$ to~\cite{Roubicek-book}. The $\Hkdiv[][f]{k}$-spaces are defined, e.g., in~\cite{brezzi-book} by \\
       $\Hkdiv[][f]{k}:=\cbrac{\Test\in \Hk{k}{}:~\nabla\cdot\Test\in \Hk{k}{}, \Test\cdot\nu = f \text{ on } \Gamma}$.%
       %Especially, we set $\Hkdiv[][0]{1}:=\Hkdiv[][0]{1}$.%
       \label{Not:spaces}
 \item \textbf{Products: } We denote by $\scp{\cdot}{\cdot}_H$ the inner product on a Hilbert space~$H$ and by $\dualp{\cdot}{\cdot}_{V^\ast\times V}$, the dual pairing between
       a Banach space~$V$ and its dual space~$V^\ast$. On $\setR^n$, we just write $\vecv\cdot\vecu:=\scp{\vecv}{\vecu}_{\setR^n}$ and on $\Lp{2}{}$,
       we just denote $\scp{\cdot}{\cdot}_{\Omega}:=\scp{\cdot}{\cdot}_{\Lp{2}{}}$. In particular the dual pairing between $\Hk{1}{}$ and its dual~$\Hk{1}{}^\ast$,
       we abbreviate by $\dualp{\cdot}{\cdot}_{1,\Omega}:=\dualp{\cdot}{\cdot}_{\Hk{1}{}^\ast\times\Hk{1}{}}$.%
       \label{Not:Prod}
\setcounter{countNotation}{\value{enumi}}
\end{enumerate}
In order to successfully examine the above model, we introduce the following assumptions
%
%=============================
% LIST: Assumptions
%=============================
\begin{enumerate}[label=({A}\arabic*), ref=({A}\arabic*), itemsep=0.0mm,  start=\value{countAssumption}+1]
 \item \textbf{Geometry: } Let $n\in\{2,3\}$ and $\Omega\subset\setR^n$ be a bounded Lipschitz domain, i.e. $\Gamma\in C^{0,1}$.%
       \label{Assump:Geom}
 \item \textbf{Initial data: } The initial data~$\clstart$ are non negative and bounded, i.e., \\
       $0\leq\clstart(x)\leq M_0$ for a.e. $x\in\Omega$ for some $M_0 \in\setR_+$.%
       \label{Assump:InitData}
 \item \textbf{Ellipticity: } The diffusivity tensor~$\Dl$ and the permeability tensor~$\permeabH$ satisfy \\
       $\Dl\xi\cdot\xi>\alpha_D\abs{\xi}^2$ and $\permeabH^{-1}\xi\cdot\xi>\alpha_K\abs{\xi}^2$ for all $\xi\in\setR^n$, \\
       $\Dl\xi\cdot\eta<C_D\abs{\xi}\abs{\eta}$ and $\permeabH^{-1}\xi\cdot\eta<C_K\abs{\xi}\abs{\eta}$ for all $\xi,\eta\in\setR^n$.%
       \label{Assump:Ellip}
 \item \textbf{Coefficients: } The porosity~$\theta$, the dynamic viscosity~$\mu$, and the electric permittivity~$\epsilon$ are positive constants.%
       \label{Assump:Coeff}
 \item \textbf{Reaction rates: } The reaction rate functions $\Rl:\setR^2\rightarrow\setR$ are global Lipschitz continuous functions, i.e., $\Rl\in \Ck[\setR^2]{0,1}{}$
       with Lipschitz constant~$C_{\Rl}$. Furthermore, we assume $\Rl(\vecnull)=0$ and $\Rl(\vecv) \geq 0$ for all $\vecv\in\setR^2$
       with $v_l\leq 0$. This means, in case a chemical species vanishes, it can only be produced.%
       \label{Assump:Reaction}
 \item \textbf{Boundary data: } We assume $\gl\in\Lp[\OmegaT]{\infty}{}$, $l=1,2$, $\sigma\in\Lp[\OmegaT]{\infty}{}$, and $f\in\Lp[\OmegaT]{\infty}{}$.
       Furthermore, we suppose that functions $\vecf,\vecsigma\in\Lp[I]{\infty}{;\Wkp{1}{\infty}{}}$ with $\vecsigma\cdot\vecnu=\sigma$ and $\vecf\cdot\vecnu=f$ exist.%
       \label{Assump:BoundData}
 \item \textbf{Background charge density: } We assume $\rho_b\in \Lp[\OmegaT]{\infty}{}$ for the background charge density.%
       \label{Assump:back-charge}
\setcounter{countAssumption}{\value{enumi}}
\end{enumerate}
Equipped with the just introduced notation, we now define the weak formulation of the model.
%
%=============================
%=============================
% DEFINITION: Weak Solution
%=============================
%=============================
\begin{definition}[Weak solution]\label{def:weaksolution}
The vector $\brac{\fieldEL,\potEL,\fieldF,\press,\Cl}\in \setR^{4+2n}$ is a weak solution of the \dpnp\ \eqref{eq:Model-1a}--\eqref{eq:Model-1j},  if and only if
\begin{subequations}
\begin{enumerate}[label=(\roman*), ref=(\roman*), itemsep=-1.5mm]
 \item $(\fieldEL,\potEL)\in \spaceSEF\times\spaceSP$ solves for all $(\Test,\test) \in \spaceTF\times\spaceTP$
      \begin{align}
        \scp{\permitEL^{-1}\fieldEL}{\Test}_\Omega &= \scp{\potEL}{\grad\cdot\Test}_\Omega  \label{eq:gaussWeak}\\
        \scp{\grad\cdot\fieldEL}{\test}_{\Omega} &= \scp{\rho_b + \chargeEL}{\test}_\Omega \label{eq:gaussWeakDIV}, \\
          \text{with the free} & \text{ charge density } \chargeEL \text{ given by }\chargeEL = \theta(\chargeELlongTwo)~ \nonumber.
      \end{align}
 \item $(\fieldF,\press)\in \spaceSF\times\spaceSP$ solves for all $(\Test,\test) \in \spaceTF\times\spaceTP$
       \begin{align}
           \scp{\permeabH^{-1}\fieldF}{\Test}_\Omega &= \scp{\mu^{-1}\press}{\grad\cdot\Test}_\Omega + \scp{\coeffForceEL\forceEL}{\Test}_\Omega, \label{eq:darcyWeak}\\
               \scp{\grad\cdot\fieldF}{\test}_\Omega &= 0~.    \label{eq:darcyWeakDIV}
       \end{align}
 \item $\cl \in \fspace{L^\infty}{I}{;\Lp{2}{}} \cap \spaceSC \cap \Lp[\OmegaT]{\infty}{}$ solves for all $\test\in\spaceTC$ and for $l=1,2$
       \begin{align}\label{eq:transportWeak}
        &\dualp{\theta\dert\cl}{\test}_{1,\Omega} + \scp{\Dl\grad\cl}{\grad\test}_\Omega - \scp{ {\cl}[\fieldF + \coeffEL\fieldEL] }{\grad\test}_\Omega \nonumber\\[2.0mm]
        &  = \scp{\theta\Rl(\Cl)}{\test}_\Omega + \scp{\gl}{\test}_\Gamma~,
       \end{align}
       and $\cl$ take its initial values in the sense that
       \begin{align*}
        \lim\limits_{t\searrow0} \scp{\cl(t) - \clstart}{\test}_\Omega ~=~0 \qquad \text{for all } \test\in\Lp{2}{}~.\\[-12mm] \nonumber
      \end{align*}
      \hfill$\square$
\end{enumerate}
\end{subequations}
\end{definition}
\begin{remark}
\textnormal{
We note that equation \eqref{eq:transportWeak} is not well-defined without having $\cl \in\Lp[\OmegaT]{\infty}{}$.
This is due to the fact, that an embedding of the type $\Hkdiv{1}{}\hookrightarrow\Lp{p}{}$, for some $p>2$, does not hold true.
Thus, we have $\fieldEL,\fieldF\in\Lp{2}{}$ at the best and for the existence of the convection integral and the electric~drift integral in \eqref{eq:transportWeak},
we need the estimate
\begin{align*}
\scp{{\cl}[\fieldF + \coeffEL\fieldEL]}{\grad\test}_\Omega ~\leq~ \norm{\cl}{\Lp{\infty}{}} \norm{\fieldF + \coeffEL\fieldEL}{\Lp{2}{}} \norm{\grad\test}{\Lp{2}{}}~.
\end{align*}
This shows that $\cl\in\Lp[\OmegaT]{\infty}{}$ is mandatory for a well-defined weak formulation. Consequently, we have to include $\Lp[\OmegaT]{\infty}{}$
in the solution space for $\cl$.%
}\hfill$\square$
\end{remark}
\begin{remark}
\textnormal{
In equations \eqref{eq:gaussWeak} and \eqref{eq:darcyWeak}, the test function space differs from the solution space and the solutions $\fieldEL$ and $\fieldF$ are not admissible test functions.
However, \eqref{eq:Model-1d}, \eqref{eq:Model-1g}, and \ref{Assump:BoundData} ensure that $\fieldEL-\vecsigma$ and $\fieldF-\vecf$ are admissible test functions.%
}\hfill$\square$
\end{remark}
%
%
%%%%%%%%%%%%%%%%%%%%%%%%%%%%%
%%%%%%%%%%%%%%%%%%%%%%%%%%%%%
%  SECTION uniqueness       %
%%%%%%%%%%%%%%%%%%%%%%%%%%%%%
%%%%%%%%%%%%%%%%%%%%%%%%%%%%%
\section{Uniqueness}\label{sec:Uniqueness}
In this section, we show that the solutions of the investigated \dpnp\ are unique.
%
%
%=========================
%=========================
%% THEOREM uniqueness
%=========================
%=========================
\begin{thm}[Uniqueness]\label{thm:unique}
Let \ref{Assump:Geom}--\ref{Assump:back-charge} be valid and let $\brac{\fieldEL,\potEL,\fieldF,\press,\Cl}\in \setR^{4+2n}$ be a weak solution of \eqref{eq:Model-1a}--\eqref{eq:Model-1j} according to \cref{def:weaksolution}.
Then, $\brac{\fieldEL,\potEL,\fieldF,\press,\Cl}$ is unique.
\end{thm}
\begin{proof}
Let us assume that $\brac{\fieldEL_i,\potEL_i,\fieldF_i,\press_i,\Cl_i}$, $i=1,2$, are two solutions of \eqref{eq:Model-1a}--\eqref{eq:Model-1j} to identical data.
Furthermore, we denote the difference between these two solution by
\begin{align*}
 \brac{\fieldEL_{12},\potEL_{12},\fieldF_{12},\press_{12},\Cl_{12}}:=\brac{\fieldEL_1-\fieldEL_2,\potEL_1-\potEL_2,\fieldF_1-\fieldF_2,\press_1-\press_2,\Cl_1-\Cl_2}
\end{align*}
By subtracting the equations for the respective solutions, we obtain the error equations\\[2.0mm]
%===================
%% error equations
%==================
\begin{subequations}
\underline{Gauss's law:}\vspace{-2mm}
\begin{align}
 \scp{\permitEL^{-1}\fieldEL_{12}}{\Test}_\Omega &= \scp{\potEL_{12}}{\grad\cdot\Test}_\Omega, \label{eq:gaussWeak-error} \\
   \scp{\grad\cdot\fieldEL_{12}}{\test}_{\Omega} &= \scp{\chargeEL[,12]}{\test}_\Omega =\theta\sum_l\scp{\zl\cl[l,12] }{\test}_\Omega~. \label{eq:gaussWeakDIV-error}
\end{align}
\underline{Darcy's law:}\vspace{-2mm}
\begin{align}
 \scp{\permeabH^{-1}\fieldF_{12}}{\Test}_\Omega &= \scp{\mu^{-1}\press_{12}}{\grad\cdot\Test}_\Omega + \theta\mu^{-1}\sum_l\scp{\zl\cl[l,1]\permitEL^{-1}\fieldEL_1-\zl\cl[l,2]\permitEL^{-1}\fieldEL_2}{\Test}_\Omega \label{eq:darcyWeak-error}\\
     \scp{\grad\cdot\fieldF_{12}}{\test}_\Omega &= 0~. \label{eq:darcyWeakDIV-error}
\end{align}
\underline{Nernst--Planck equations:}\vspace{-2mm}
\begin{align}
   & \dualp{\theta\dert\cl[l,12]}{\test}_{1,\Omega} + \scp{\Dl\grad\cl[l,12]}{\grad\test}_\Omega
     - \scp{ {\cl[l,1]}[\fieldF_1 + \coeffEL\fieldEL_1] }{\grad\test}_\Omega \nonumber\\[2.0mm]
   & + \scp{ {\cl[l,2]}[\fieldF_2 + \coeffEL\fieldEL_2]  }{\grad\test}_\Omega = \theta \scp{\Rl(\Cl_1)-\Rl(\Cl_2)}{\test}_\Omega ~. \label{eq:transportWeak-error}
\end{align}
\end{subequations}
%
%========================
%% Nernst-Planck - error
%========================
We now show by contradiction that $\cl[l,12]\equiv0$ for $l=1,2$. To this end, we assume that
\begin{align}\label{eq:unique-contradiction}
 \sum_l \norm{\cl[l,12]}{\Lp{2}{}}^2 > 0 ~~\Equivalent~~ \exists~\kappa>0 \text{ such that } \sum_l \norm{\cl[l,12]}{\Lp{2}{}}^2 \geq \kappa~.
\end{align}
Next, we test equation~\eqref{eq:transportWeak} with $\test=\cl[l,12]$ and we sum over $l=1,2$. Thereby, we come for the time integral
and the diffusion integral with \ref{Assump:Ellip}  to
\begin{align*}
      \sum_l \dualp{\theta\dert\cl[l,12]}{\cl[l,12]}_{1,\Omega} + \scp{\Dl\grad\cl[l,12]}{\grad\cl[l,12]}_\Omega
 \geq \frac{\theta}{2}\derr \sum_l\norm{\cl[l,12]}{\Lp{2}{}}^2 + \alpha_D \sum_l\norm{\grad\cl[l,12]}{\Lp{2}{}}^2.
\end{align*}
For the reaction integral, we obtain with \ref{Assump:Reaction}
\begin{align*}
      \theta\sum_l \scp{\Rl(\Cl_1)-\Rl(\Cl_2)}{\cl[l,12]}_\Omega \leq \theta\sum_l C_{\Rl}\norm{ \abs{\Cl_{12}} \cl[l,12]}{\Lp{1}{}}
 \leq \theta\max_l C_{\Rl} \sum_l \norm{\cl[l,12]}{\Lp{2}{}}^2.
\end{align*}
The convection integral and the electric~drift integral, we transform to
\begin{align*}
 &     -\sum_l \scp{ {\cl[l,1]}[\fieldF_1 + \coeffEL\fieldEL_1] - {\cl[l,2]}[\fieldF_2 + \coeffEL\fieldEL_2] }{\grad\cl[l,12]}_\Omega \\
 &=    -\sum_l \scp{ {\cl[l,1]}[\fieldF_{12} + \coeffEL\fieldEL_{12}]}{\grad\cl[l,12]}_\Omega \\
 &~~~~ -\sum_l \scp{ {\cl[l,12]}[\fieldF_2 + \coeffEL\fieldEL_2] }{\grad\cl[l,12]}_\Omega
  ~~=: A.1 +  A.2~,
\end{align*}
and $A.2$, we estimate with \eqref{eq:unique-contradiction} and  \Young\ with a free parameter $\delta>0$, cf. \cite{GilbargTrudinger-book}, by
\begin{align*}
 A.2 &\leq \delta\sum_l \norm{\grad\cl[l,12]}{\Lp{2}{}}^2
           +\frac{\kappa}{\kappa}\frac{1}{2\delta}\norm{\fieldF_2 + \coeffEL\fieldEL_2}{\Lp{2}{}}^2\sum_l\norm{\cl[l,12]}{\Lp{\infty}{}}^2\\
%      &\leq \delta\sum_l \norm{\grad\cl[l,12]}{\Lp{2}{}}^2
%            + \frac{1}{2\kappa\delta}\norm{\fieldF_2 + \coeffEL\fieldEL_2}{\Lp{2}{}}^2\norm{\Cl_{12}}{\Lp{\infty}{}}^2\sum_l\norm{\cl[l,12]}{\Lp{2}{}}^2 \\
     &\leq \delta\sum_l \norm{\grad\cl[l,12]}{\Lp{2}{}}^2
           + \frac{1}{2\kappa\delta}\sum_i\sqbrac{\norm{\fieldF_i + \coeffEL\fieldEL_i}{\Lp{2}{}}^2\norm{\Cl_i}{\Lp{\infty}{}}^2}\sum_l\norm{\cl[l,12]}{\Lp{2}{}}^2.
\end{align*}
Analogously, we come for $A.1$ to
\begin{align*}
 A.1 &\leq \delta\sum_l \norm{\grad\cl[l,12]}{\Lp{2}{}}^2  + \frac{2}{\kappa\delta}\sum_i\sqbrac{\norm{\fieldF_i + \coeffEL\fieldEL_i}{\Lp{2}{}}^2\norm{\Cl_i}{\Lp{\infty}{}}^2}\sum_l\norm{\cl[l,12]}{\Lp{2}{}}^2.
\end{align*}
Altogether, we arrive with a proper choice of a free parameter~$\delta>0$ at
\begin{align*}
 &     \frac{\theta}{2}\derr \sum_l\norm{\cl[l,12]}{\Lp{2}{}}^2 + \frac{\alpha_D}{2} \sum_l\norm{\grad\cl[l,12]}{\Lp{2}{}}^2 \\
 &\leq \brac{\theta\max_l C_{\Rl}+\frac{4}{\kappa\alpha_D}\sum_i\sqbrac{\norm{\fieldF_i + \coeffEL\fieldEL_i}{\Lp{2}{}}^2\norm{\Cl_i}{\Lp{\infty}{}}^2}}\sum_l\norm{\cl[l,12]}{\Lp{2}{}}^2.
\end{align*}
We note that the initial values vanish due to $\cl[l,12]=\clstart-\clstart=0$. Thus, applying \Grownwall, cf. \cite{Evans-book}, yields
\begin{align*}
 \sum_l\norm{\cl[l,12](t)}{\Lp{2}{}}^2 \leq C \sum_l\norm{\cl[l,12](0)}{\Lp{2}{}}^2 =0 \quad \text{for a.e. } t\in[0,\timeEnd] ~.
\end{align*}
This is a contradiction to assumption \eqref{eq:unique-contradiction}. Hence, we have proven $\cl[l,12]\equiv0$.
\par
%
%===================
%% Gauss - error
%==================
We proceed by testing equation~\eqref{eq:gaussWeakDIV-error} with $\test=\grad\cdot\fieldEL_{12}$. Thereby, we come with \Young\ and  $\cl[l,12]\equiv0$ directly to
\begin{align*}
  \norm{\grad\cdot\fieldEL_{12}}{\Lp{2}{}}^2 \leq \theta\max_l\abs{\zl} \sum_l\norm{\cl[l,12]}{\Lp{2}{}}^2=0~.
\end{align*}
Next, we test equation~\eqref{eq:gaussWeakDIV-error} with $\test=\potEL_{12}$ and equation~\eqref{eq:gaussWeak-error} with $\Test=\fieldEL_{12}$.
By adding these equations, we get with \ref{Assump:Ellip} and $\cl[l,12]\equiv0$
\begin{align*}
  \permitELscalar\alpha_D\norm{\fieldEL_{12}}{\Lp{2}{}}^2 \leq \scp{\permitEL^{-1}\fieldEL_{12}}{\fieldEL_{12}}_\Omega = \theta \sum_l\scp{\zl\cl[l,12] }{\potEL_{12}}_\Omega = 0~.
\end{align*}
We now test equation~\eqref{eq:gaussWeak-error} with $\Test\in\spaceTF$, for which we assume that $\grad\cdot\Test=\potEL_{12}$
and $\norm{\Test}{\Hkdiv{1}{}}\leq C\norm{\potEL_{12}}{\spaceTP}$ holds, cf. \cite[Chapter~7.2]{Quarteroni-book}.
This gives with \ref{Assump:Ellip} and \Young\
\begin{align*}
   \norm{\potEL_{12}}{\Lp{2}{}}^2 \leq \frac{2}{\permitELscalar\alpha_D}\norm{\fieldEL_{12}}{\Lp{2}{}}^2 = 0~.
\end{align*}
Hence, we have proven $\norm{\potEL_{12}}{\Lp{2}{}}^2 + \norm{\fieldEL_{12}}{\Hkdiv{1}{}}^2 \leq 0$, which means $\potEL_{12}\equiv0, \fieldEL_{12}\equiv\vecnull, \grad\cdot\fieldEL_{12}\equiv0$.
\par
%
%===================
%% Darcy - error
%==================
Analogously, we test equation~\eqref{eq:darcyWeakDIV-error} with $\test=\grad\cdot\fieldF_{12}$. This shows $\norm{\grad\cdot\fieldF_{12}}{\Lp{2}{}}^2 =0$.
Next, we test equation~\eqref{eq:darcyWeakDIV-error} with $\test=\press_{12}$ and equation~\eqref{eq:darcyWeak-error} with $\Test=\fieldF_{12}$.
Then, we add these equations and obtain with \ref{Assump:Ellip}, $\cl[l,12]\equiv0$, and $\fieldEL_{12}\equiv\vecnull$
\begin{align*}
         \alpha_K \norm{\fieldF_{12}}{\Lp{2}{}}^2
  &\leq  \theta\mu^{-1}\sum_l\scp[]{\zl\cl[l,12]\permitEL^{-1}\fieldEL_1 + \zl\cl[l,2]\permitEL^{-1}\fieldEL_{12}}{\fieldF_{12}} = 0~.
%  %% LONG VERSION
%          \alpha_K \norm{\fieldF_{12}}{\Lp{2}{}}^2
% %   &\leq %\scp[]{\permeabH^{-1}\fieldF_{12}}{\fieldF_{12}}_\Omega =
% %         \theta\mu^{-1}\sum_l\scp[]{\zl\cl[l,1]\permitEL^{-1}\fieldEL_1 - \zl\cl[l,2]\permitEL^{-1}\fieldEL_2}{\fieldF_{12}}_\Omega \\
%   &=    \theta\mu^{-1}\sum_l\scp[]{\zl\cl[l,12]\permitEL^{-1}\fieldEL_1 + \zl\cl[l,2]\permitEL^{-1}\fieldEL_{12}}{\fieldF_{12}} = 0~.
\end{align*}
By testing equation~\eqref{eq:darcyWeak-error} with $\Test\in\spaceTF$, for which we assume according to \cite[Chapter~7.2]{Quarteroni-book} that
$\grad\cdot\Test=\press_{12}$ and $\norm{\Test}{\Hkdiv{1}{}}\leq C\norm{\press_{12}}{\spaceTP}$ holds, we come with \ref{Assump:Ellip} and \Young\ to
\begin{align*}
  \norm{\press_{12}}{\Lp{2}{}}^2 &\leq \delta C\norm{\press_{12}}{\Lp{2}{}}^2 + \frac{C_k^2}{2\delta}\norm{\fieldF_{12}}{\Lp{2}{}}^2
                                        + \frac{\theta}{2\delta\mu}\norm{\chargeEL[,1]\permitEL^{-1}\fieldEL_1 - \chargeEL[,2]\permitEL^{-1}\fieldEL_2}{\Lp{2}{}}^2 \\
                                 &=:I.1+I.2+I.3~.
\end{align*}
We already know $I.2=0$ and $I.3$ is estimated with $\cl[l,12]\equiv0$ and $\fieldEL_{12}\equiv\vecnull$ by
\begin{align*}
 I.3 &= \norm{\chargeEL[,12]\permitEL^{-1}\fieldEL_1 + \chargeEL[,2]\permitEL^{-1}\fieldEL_{12}}{\Lp{2}{}}^2 = 0~.
\end{align*}
Thus, we have $\norm{\press_{12}}{\Lp{2}{}}=0$ and we have proven $\norm{\press_{12}}{\Lp{2}{}}^2 + \norm{\fieldF_{12}}{\Hkdiv{1}{}}^2 \leq 0$, which means
$\press_{12}\equiv0, \fieldF_{12}\equiv\vecnull, \grad\cdot\fieldF_{12}\equiv0$.
\end{proof}
%
%
%%%%%%%%%%%%%%%%%%%%%%%%%%%%%%%%%%
%%%%%%%%%%%%%%%%%%%%%%%%%%%%%%%%%%
%  SECTION fixed point operator  %
%%%%%%%%%%%%%%%%%%%%%%%%%%%%%%%%%%
%%%%%%%%%%%%%%%%%%%%%%%%%%%%%%%%%%
\section{Fixed Point Operator} \label{sec:fpi-operator}
In the next sections, we prove the existence of solutions of the \dpnp\ by applying a fixed point approach. The idea behind this method of proof can be roughly summarized as follows:
\par
Firstly, linearize the nonlinear system with a suitable linearization method. For that purpose, often well-known and widely used numerical linearization schemes are used.
Concerning \dpnp{s}, the most famous linearization scheme for numerical computations is the so-called Gummel iteration, cf. \cite{Gummel1964}.
\par
Secondly, reformulate the linearized system by means of an abstract operator. This operator is exactly constructed such that the images of this operator are
the solutions of the linearized system. Furthermore, the construction must by carried out in such way, that the solutions of the nonlinear system are exactly
the fixed points of this operator. Hence, the existence of solutions of the nonlinear system is equivalent to the existence of fixed points of the constructed operator.
\par
Thirdly, it remains to prove that the operator satisfies the assumptions of a fixed point theorem, which allows to conclude that a fixed point exists.
This is the reason why the most part of the subsequent proof consists in verifying the assumptions of the fixed point~\cref{thm:fpi-theorem}.
\medskip
\par
%===================================================
% LINEARIZATION: Gummel type linearization approach
%===================================================
We now linearize the \dpnp\ by the following Gummel-type approach, which is sketched as follows.
\begin{enumerate}[label=(L.\arabic*), ref=(L.\arabic*),topsep=2.0mm, itemsep=-1.0mm, start=0]
 \item We replace the free charge density~$\chargeEL$ by some given approximation~$\chargeELold$.%
       \label{linStep:start}
 \item Thereby, we decouple Gauss's law from the remaining \dpnp, as for a given $\chargeELold$, we obtain a solution $(\fieldEL,\potEL)$ of  Gauss's law independently of the
       remaining solution vector~$(\fieldF,\press,\Cl)$. In the following \cref{def:FixOp}, this is formulated by means of the operator~$\fixOp_1$.%
       \label{linStep:gauss}
 \item Next, we proceed by solving Darcy's law. From \ref{linStep:start} and \ref{linStep:gauss} we know that we can take $(\chargeELold,\fieldEL,\potEL)$ as a given input.
       Hence, we obtain a solution~$(\fieldF,\press)$ independently of the remaining solution vector~$\Cl$. In the following \cref{def:FixOp}, this is formulated by means of
       the solution operator~$\fixOp_2$.%
       \label{linStep:darcy}
 \item Finally, we know from \ref{linStep:start}--\ref{linStep:darcy}, that we can treat $(\chargeELold,\fieldEL,\potEL,\fieldF,\press)$ as a given input for the equations for $\cl$,
       which gives immediately the remaining solution~$\Cl$. This is formulated in the following \cref{def:FixOp} by means of the solution operator~$\fixOp_3$.%
       \label{linStep:np}
\end{enumerate}
This Gummel-type linearization approach is rigorously formulated in the next definition.
%
%
%==================================
%==================================
% DEFINITION: fixed point operator
%==================================
%==================================
\begin{definition}[Fixed point operator]\label{def:FixOp}
Let $K\subset X$ be a subset of the Banach space $X$, which is given by $X:= \sqbrac{\fspace{L^\infty}{I}{;\Lp{2}{}} \cap \spaceSC \cap \Lp[\OmegaT]{\infty}{}}^2$.
We introduce the fixed point operator~$\fixOp$ by
\begin{align*}
\fixOp:=\fixOp_3\circ\fixOp_2\circ\fixOp_1: K \subset X \rightarrow X~.
\end{align*}
Herein, the suboperator $\fixOp_1$ is defined by
\begin{subequations}
\begin{align}
 \fixOp_1: &\begin{cases} \fixSet \rightarrow X\times\spaceSEF\times\spaceSP ~=:Y \\
                           \Clold ~\mapsto (\Clold,\fieldEL,\potEL), ~\text{ with } (\fieldEL,\potEL) \text{ solving for all } (\Test,\test) \in \spaceTF\times\spaceTP
              \end{cases} \nonumber\\
            & \hspace{27.0mm} \scp{\permitEL^{-1}\fieldEL}{\Test}_\Omega = \scp{\potEL}{\grad\cdot\Test}_\Omega , \label{eq:gaussWeak-fpi}    \\
            & \hspace{27.0mm} \scp{\grad\cdot\fieldEL}{\test}_{\Omega} = \scp{\rho_b + \chargeELold}{\test}_\Omega , \label{eq:gaussWeakDIV-fpi} \\
            & \hspace{27.0mm} \text{with the free charge density } \chargeELold \text{ given by }\chargeELold = \theta(\chargeELlongTwoold)~. \nonumber
\end{align}
Furthermore, the suboperator $\fixOp_2$ is defined by
\begin{align}
 \fixOp_2: &\begin{cases} \qquad Y ~~~~\rightarrow Y\times \spaceSF\times\spaceSP~=:Z \\
                          (\Clold,\fieldEL,\potEL) \mapsto (\Clold,\fieldEL,\potEL,\fieldF,\press), ~\text{ with } (\fieldF,\press) \text{ solving for all } (\Test,\test) \in \spaceTF\times\spaceTP
              \end{cases} \nonumber\\
            & \hspace{27.0mm} \scp{\permeabH^{-1}\fieldF}{\Test}_\Omega = \scp{\mu^{-1}\press}{\grad\cdot\Test}_\Omega + \scp{\coeffForceEL\forceELold}{\Test}_\Omega, \label{eq:darcyWeak-fpi}   \\
            & \hspace{27.0mm} \scp{\grad\cdot\fieldF}{\test}_\Omega = 0 .   \label{eq:darcyWeakDIV-fpi}
\end{align}
Finally, the suboperator $\fixOp_3$ is defined by
\begin{align}
 \fixOp_3: &\begin{cases} \qquad~~~ Z \qquad~~\rightarrow X \\
                          (\Clold,\fieldEL,\potEL,\fieldF,\press) \mapsto \Cl=(\cl[1],\cl[2]), ~\text{ with } \cl \text{ solving for all } \test\in\spaceTC \text{ and } l=1,2
              \end{cases} \nonumber\\
            & \hspace{27.0mm} \dualp{\theta\dert\cl}{\test}_{1,\Omega} + \scp{\Dl\grad\cl}{\grad\test}_\Omega - \scp{ {\cl}[\fieldF + \coeffEL\fieldEL] }{\grad\test}_\Omega \nonumber\\[2.0mm]
            & \hspace{27.0mm} = \scp{\theta\Rl(\Cl)}{\test}_\Omega + \scp{\gl}{\test}_\Gamma , \label{eq:transportWeak-fpi} \\
            & \hspace{27.0mm} \text{and } \cl \text{ take its initial values in the sense that } \nonumber\\
            & \hspace{27.0mm} \lim\limits_{t\searrow0} \scp{\cl(t) - \clstart}{\test}_\Omega ~=~0 \qquad \text{for all } \test\in\Lp{2}{}~. \nonumber\\[-12mm] \nonumber
\end{align}
\end{subequations}
\hfill$\square$
\end{definition}
\begin{remark}
\textnormal{
We note, that the fixed point operator~$\fixOp$ is solely a function of~$\Clold$. For this reason, a fixed point~$\Cl$ of $\fixOp$ is only a partial solution
in the sense of \cref{def:weaksolution}, as $\Cl$ only solves the equations~\eqref{eq:transportWeak}.
However, the suboperators~$\fixOp_1$, $\fixOp_2$, and $\fixOp_3$ contain the necessary information about the remaining partial solutions $(\potEL,\fieldEL)$ and $(\press,\fieldF)$.
Furthermore, in case a fixed point $\Cl=\fixOp(\Cl)$ exists, these supoperators ensure the existence of the partial solutions $(\potEL,\fieldEL)$ and $(\press,\fieldF)$
such that this yields the existence of a solution $\brac{\fieldEL,\potEL,\fieldF,\press,\Cl}\in \setR^{4+2n}$ according to \cref{def:weaksolution}.
%Hence, it is admissible to focus on the concentrations~$c^\pm$.
}\hfill$\square$
\end{remark}
%
%
%============================================
%============================================
% LEMMA: well-defined fixed point operator
%============================================
%============================================
\begin{lemma}[well-definedness]\label{lemma:wellDef-fixOp}
Let \ref{Assump:Geom}--\ref{Assump:back-charge} be valid. Then, $\fixOp: \fixSet\subset X \rightarrow X~$ defined in \cref{def:FixOp}, is well-defined.
\end{lemma}
\begin{proof}
 $\fixOp_1$ is well-defined in the first component, since $\fixOp_1$ is the identity in the first component. As to the components $(\fieldEL,\potEL)$,
 we know that for all $\Cl\in K$ unique solutions $(\fieldEL,\potEL)\in\spaceSEF\times\spaceSP$ of \eqref{eq:gaussWeak-fpi} and \eqref{eq:gaussWeakDIV-fpi} exist.
 This follows from \cite[Theorem 7.4.1]{Quarteroni-book}.
 However, $\potEL$ is only determined up to a constant. Imposing, e.g.,  a zero mean value constraint%
 \footnote{\label{foot1}The mean value of a function $f\in\Lp{1}{}$ is defined by $\frac{1}{\abs{\Omega}}\Intdx{ f~}$.},
 leads to uniqueness of $\potEL$. Furthermore, we note that the time variable $t$ plays only the role of a parameter in the equations for $(\fieldEL,\potEL)$.
 This leads to uniform results with respect to $t$. Hence, $\fixOp_1$ is well-defined.
 \par
 $\fixOp_2$ is the identity in the first three components. For the last two components $(\fieldF,\press)$,
 we know that for all $(\Cl,\fieldEL)\in K\times\spaceSEF$ unique solutions of \eqref{eq:darcyWeak-fpi} and \eqref{eq:darcyWeakDIV-fpi} exist.
 This follows again from \cite[Theorem 7.4.1]{Quarteroni-book}.
 Likewise, $\press$ is only determined up to a constant and we obtain uniqueness by imposing, e.g.,  a zero mean value constraint.
 The existence is uniform in time, as $t$ plays just the role of a parameter. Thus, $\fixOp_2$ is well-defined.
 \par
 Applying Rothe's method, cf. \cite[Chapter 8.2]{Roubicek-book}, \cite{rektorys-book}, together with the regularities of $\fieldEL,\potEL$, $\fieldF,\press$ (according to \cref{thm:aprioriBounds})
 guarantees with \cref{lemma:nonnegative} and \cref{lemma:bounded} the existence of unique weak solutions $\Cl\in X$ of equations~\eqref{eq:transportWeak-fpi}.
 Thus, $\fixOp_3$ is well defined.
\end{proof}
%
%
%============================================
%============================================
% LEMMA: regularity for electric field
%============================================
%============================================
\begin{lemma}[Regularity for Gauss's law]\label{lemma:regularity-gausslaw}
Let \ref{Assump:Geom}--\ref{Assump:back-charge} be valid and let $(\fieldEL,\potEL,\fieldF,\press,\Cl)\in\setR^{4+2n}$ be a solution of \eqref{eq:Model-1a}--\eqref{eq:Model-1j} according to \cref{def:FixOp}.
Then, for the partial solution~$(\fieldEL,\potEL)$, we have
\begin{align*}
 \potEL\in\Lp[I]{\infty}{;\Hk{2}{}/\setR} \qquad \text{ and } \qquad \fieldEL\in\Lp[I]{\infty}{;\Hk{1}{}}~.
\end{align*}
\end{lemma}
\begin{proof}
 We recall from \cite{Grisvard-book,Quarteroni-book}, that the equation %\cite[Chapter~2]{Grisvard-book}, \cite[Chapter~6.1]{Quarteroni-book}
 \begin{align*}
  \scp{\permitEL\grad\psi}{\grad\test}_{\Omega} = \scp{\rho_b+\chargeELold}{\test}_\Omega + \scp{\sigma}{\test}_\Gamma \qquad \text{ for all } \test\in\Hk{1}{}.
 \end{align*}
 possesses a unique solution~$\psi\in\Hk{2}{}$, if we impose a zero mean value constraint\footnote{See footnote \ref{foot1}}.
%  More precisely, we obtain a solution $\psi\in\Lp[I]{\infty}{;\Hk{1}{}}$ and by standard elliptic regularity theory we get $\psi\in\Lp[I]{\infty}{;\Hk{2}{}}$.
 As the time variable~$t$ plays only the role of a parameter in the equation for $\psi$, we obtain all results uniformly in time. This yields $\psi\in\Lp[I]{\infty}{;\Hk{2}{}}$.
 Hence, by defining
 \begin{align*}
  \potEL:=\psi\in\Lp[I]{\infty}{;\Hk{2}{}} \qquad\text{ and }\qquad \fieldEL:=\permitEL\grad\psi\in\Lp[I]{\infty}{;\Hk{1}{}},
 \end{align*}
 we obtain a solution of equations~\eqref{eq:gaussWeak-fpi} and \eqref{eq:gaussWeakDIV-fpi}. Finally, we already know from the proof of \cref{lemma:wellDef-fixOp},
 that the above constructed solution~$(\fieldEL,\potEL)$ is the unique solution of equations~\eqref{eq:gaussWeak-fpi} and \eqref{eq:gaussWeakDIV-fpi}.
\end{proof}
%
%
%%%%%%%%%%%%%%%%%%%%%%%%%%%%%%%%%%%%%%%%%%%%%%%%%%%%%%%
%%%%%%%%%%%%%%%%%%%%%%%%%%%%%%%%%%%%%%%%%%%%%%%%%%%%%%%
%  SECTION existence for two-component electrolytes   %
%%%%%%%%%%%%%%%%%%%%%%%%%%%%%%%%%%%%%%%%%%%%%%%%%%%%%%%
%%%%%%%%%%%%%%%%%%%%%%%%%%%%%%%%%%%%%%%%%%%%%%%%%%%%%%%
\section{Global Existence of a Solution} \label{sec:Existence-twoComp}
%
%
%%%%%%%%%%%%%%%%%%%%%%%%%%%%%%%%
%  SUBSECTION A priori bounds  %
%%%%%%%%%%%%%%%%%%%%%%%%%%%%%%%%
\subsection{A priori Estimates} \label{subsec:aprioriBounds}
In this section, we show a~priori bounds for the solution vector $\brac{\fieldEL,\potEL,\fieldF,\press,\Cl}\in \setR^{4+2n}$.
We begin with some preliminary results, which we need throughout the rest of  this paper.
Henceforth, we denote by $C$ a generic constant, which may change from line to line in the calculations.
%
%
%=======================
%=======================
% LEMMA: sign condition
%=======================
%=======================
\begin{lemma}[Algebraic Inequality]\label{lemma:sign}
 Let $p\geq0$ and $a,b\in\setR$ with $a\geq0$ and $b\geq0$. Then, we have
 \begin{align*}
  (a-b)(a^p-b^p)=(b-a)(b^p-a^p) \geq 0~.
 \end{align*}
\end{lemma}
\begin{proof}
 The equality is obvious and the inequality follows by considering the cases $a\geq b$ and $a<b$.
\end{proof}
%
%
%===============================
%===============================
% LEMMA: Boundary interpolation
%===============================
%===============================
\begin{lemma}[Boundary Interpolation]\label{lemma:interpol-boundary}
 Let $u\in\Hk{1}{}$ and suppose \ref{Assump:Geom}. Then, we have
 \begin{align*}
  \norm{u}{\Lp[\Gamma]{2}{}}^2 ~\leq~ \delta \norm{\grad u}{\Lp{2}{}}^2 + 2\delta^{-1} \norm{u}{\Lp{2}{}}^2 \qquad \text{ for all } \delta\in(0,1)~.
 \end{align*}
\end{lemma}
\begin{proof}
 Let $s\in [1/2,1)$ and $\Hk{s}{}$ be a fractional Sobolev space defined in \cite[Chapter~7.35, 7.43]{Adams1-book}. According to \cite[Theorem~7.58]{Adams1-book}
and \cite[Lemma~7.16]{Adams1-book}, we have the embeddings
\begin{align*}
   \norm{u}{\Lp[\Gamma]{2}{}} \leq C \norm{u}{\Hk{s}{}} ~\leq~ C \norm{u}{\Hk{1}{}}^s \norm{u}{\Lp{2}{}}^{1-s} ~~\text{ for } s \in [1/2,1).
\end{align*}
Then, we choose $s=1/2$ and apply \Young\ with a free parameter $\delta\in(0,1)$.
\end{proof}
We now show a lower bound for the chemical species~$\cl$.
%
%
%=======================
%=======================
% LEMMA: nonnegative
%=======================
%=======================
\begin{lemma}[Non negativity]\label{lemma:nonnegative}
Let \ref{Assump:Geom}--\ref{Assump:back-charge} be valid and let $\brac{\fieldEL,\potEL,\fieldF,\press,\Cl}\in \setR^{4+2n}$
be a weak solution  of \eqref{eq:Model-1a}--\eqref{eq:Model-1j} according to  \cref{def:FixOp}. Then, we have for $l\in\{1,2\}$
\begin{align*}
\cl\tx \geq 0 \quad \text{ for a.e. } t\in[0,\timeEnd],~ \text{ a.e. } x\in\Omega~.
\end{align*}
\end{lemma}
\begin{proof}
We note that the following proof is independent of $\clold\neq\cl$ or $\clold=\cl$. For $l=1,2$, we modify equations~\eqref{eq:transportWeak-fpi} with $\cl[l,+]:=\max(\cl,0)$ to
\begin{align}\label{eq:transportWeak-modified}
  &\dualp{\theta\dert\cl}{\test}_{1,\Omega} + \scp{\Dl\grad\cl}{\grad\test}_\Omega - \scp{ {\cl[l,+]}[\fieldF + \coeffEL\fieldEL] }{\grad\test}_\Omega \nonumber\\[2.0mm]
  &  = \scp{\theta\Rl(\Cl)}{\test}_\Omega + \scp{\gl}{\test}_\Gamma~.
\end{align}
Obviously, equations~\eqref{eq:transportWeak-fpi} and \eqref{eq:transportWeak-modified} are identical for nonnegative solutions~$\cl$.
This means that nonnegative solutions~$\cl$ of \eqref{eq:transportWeak-modified} are solutions of \eqref{eq:transportWeak-fpi}.
Furthermore, by involving \cref{thm:unique}, we know that nonnegative solutions~$\cl$ of \eqref{eq:transportWeak-modified}
are the unique solutions of equations~\eqref{eq:transportWeak-fpi}. Hence, it suffices to show that \eqref{eq:transportWeak-modified} solely allows for nonnegative solutions.
\par
To this end, we test \eqref{eq:transportWeak-modified} with $\cl[l,-]:=\min(\cl,0)$. Thereby, we obtain for the time integral and the diffusion integral with \ref{Assump:Ellip}
\begin{align*}
      \sum_l \dualp{\theta\dert\cl}{\cl[l,-]}_{1,\Omega} + \scp{\Dl\grad\cl}{\grad\cl[l,-]}_\Omega
 \geq \frac{\theta}{2}\derr \sum_l\norm{\cl[l,-]}{\Lp{2}{}}^2 + \alpha_D \sum_l\norm{\grad\cl[l,-]}{\Lp{2}{}}^2.
\end{align*}
The convection integral and the electric~drift integral vanish due to $\Omega\cap\cbrac{\cl<0}\cap\cbrac{\cl>0}=\varnothing$.
For the reaction integrals and the surface integrals, we come with \ref{Assump:Reaction}, \ref{Assump:BoundData}, \cref{lemma:interpol-boundary}, and \Holder\ to
\begin{align*}
       \theta\sum_l \scp{\Rl(\Cl)}{\cl[l,-]}_\Omega + \scp{\gl}{\cl[l,-]}_\Gamma
  \leq \norm{\gl}{\Lp[\Gamma]{\infty}{}}\sum_l\norm{\cl[l,-]}{\Lp[\Gamma]{1}{}}
  \leq C\norm{\gl}{\Lp[\Gamma]{\infty}{}}\sum_l\norm{\cl[l,-]}{\Lp{2}{}} .
\end{align*}
Combining the previous estimates leads us to
\begin{align*}
      \frac{\theta}{2}\derr \sum_l\norm{\cl[l,-]}{\Lp{2}{}}^2 + \frac{\alpha_D}{2} \sum_l\norm{\grad\cl[l,-]}{\Lp{2}{}}^2
 \leq C\norm{\gl}{\Lp[\Gamma]{\infty}{}}\sum_l\norm{\cl[l,-]}{\Lp{2}{}} .
\end{align*}
It is either $\sum_l\norm{\cl[l,-]}{\Lp{2}{}}=0$ and we are done, or we have $\sum_l\norm{\cl[l,-]}{\Lp{2}{}}\geq \kappa$, for some $\kappa>0$. This gives
\begin{align*}
      \frac{\theta}{2}\derr \sum_l\norm{\cl[l,-]}{\Lp{2}{}}^2 + \frac{\alpha_D}{2} \sum_l\norm{\grad\cl[l,-]}{\Lp{2}{}}^2
 \leq \kappa^{-1}C\norm{\gl}{\Lp[\Gamma]{\infty}{}}\sum_l\norm{\cl[l,-]}{\Lp{2}{}}^2 .
\end{align*}
Applying \Grownwall\ and \ref{Assump:InitData} immediately yields $\sum_l\norm{\cl[l,-]}{\Lp{2}{}}=0$.
\end{proof}
Next, we prove energy estimates for the chemical species~$\cl$, by using the above mentioned weighted test functions, see \cref{sec:Introduction}. These energy estimates are crucial for all following results.
%
%
%============================
%============================
% LEMMA: Energy Estimates
%============================
%============================
\begin{lemma}[Energy estimates]\label{lemma:energy}
Let \ref{Assump:Geom}--\ref{Assump:back-charge} be valid and let $\brac{\fieldEL,\potEL,\fieldF,\press,\Cl}\in \setR^{4+2n}$ be a weak solution
of \eqref{eq:Model-1a}--\eqref{eq:Model-1j} according to \cref{def:FixOp}. Then, we have
\begin{align*}
  \sum_l\sqbrac{ \norm{\cl}{\Lp[I]{\infty}{;\Lp{2}{}}}+\norm{\cl}{\Lp[I]{2}{;\spaceTC}} }~\leq~ C_0
\end{align*}
Herein, the dependency of the constant is
\begin{align*}
 C_0=C_0\!\brac{\timeEnd,\max_l\abs{\zl},\norm{\gl}{\Lp[\OmegaT]{\infty}{}},\norm{f}{\Lp[\GammaT]{\infty}{}},\norm{\sigma}{\Lp[\GammaT]{\infty}{}},\norm{\rho_b}{\Lp[\OmegaT]{\infty}{}},\norm{\clstart}{\Lp{2}{}} }.
\end{align*}
\end{lemma}
\begin{proof}
For ease of readability, we split the proof into two cases.\\[2.0mm]
%=======================================================
% ENERGY - CASE 1: \cl=\clold
%=======================================================
\underline{Case 1: $\cl=\clold$}~~
In equations~\eqref{eq:transportWeak-fpi}, we choose the test functions $\varphi := \abs{\zl}\cl \in \spaceTC$ and we sum over $l=1,2$.
Thereby, we get for the time integrals and the diffusion integrals with \ref{Assump:Ellip}
\begin{align*}
 &     \sum_l \dualp{\theta\dert\cl}{\abs{\zl}\cl}_{1,\Omega} + \sum_l\scp{\Dl\grad\cl}{\grad(\abs{\zl}\cl)}_\Omega \\
 &\geq \frac{\theta}{2}\derr\sum_l \abs{\zl}\norm{\cl}{\Lp{2}{}}^2 + \alpha_D \sum_l\abs{\zl}\norm{\grad\cl}{\Lp{2}{}}^2~~.
\end{align*}
For the convection integrals, we firstly use integration by parts and we insert equation \eqref{eq:darcyWeakDIV-fpi}. Secondly, we use \Holder\ and \cref{lemma:interpol-boundary}
with a rescaled free parameter $\delta=\norm{f}{\Lp[\GammaT]{\infty}{}}^{-1}\delta$. This leads us to
\begin{align*}
 &     -\sum_l\scp{{\cl}\fieldF}{\grad(\abs{\zl}\cl)}_\Omega
  =    -\frac{1}{2} \sum_l\abs{\zl}\scp{\fieldF}{\grad\cl^2}_\Omega
  =    -\frac{1}{2} \sum_l\abs{\zl}\scp{f}{\cl^2}_\Gamma \\
 &\geq -\delta \sum_l\abs{\zl}\norm{\grad\cl}{\Lp{2}{}}^2
  ~-~  \delta^{-1} \norm{f}{\Lp[\GammaT]{\infty}{}}^2 \sum_l\abs{\zl}\norm{\cl}{\Lp{2}{}}^2 .
\end{align*}
Analogously, for the electric~drift integral we integrate by parts and we insert equation \eqref{eq:gaussWeakDIV-fpi}. This yields
\begin{align*}
 I_{el}:=&  -\frac{e}{\permitELscalar\boltz\temp} \sum_l\scp{\zl{\cl}\fieldEL}{\grad(\abs{\zl}\cl)}_\Omega
            ~= -\frac{e}{\permitELscalar\boltz\temp} \sum_l\Sign(\zl)\scp{\abs{\zl}\cl\fieldEL}{\grad(\abs{\zl}\cl)}_\Omega\\
         &=  \frac{e}{\permitELscalar\boltz\temp} \sum_l\Sign(\zl)\sqbrac{~\scp{\rho_b+\chargeELold}{(\abs{\zl}\cl)^2}_\Omega - \scp{\sigma}{(\abs{\zl}\cl)^2}_\Gamma~}.
\end{align*}
Together with \Holder\ and \cref{lemma:interpol-boundary}, we reach from this identity (with a rescaled $\delta$) at
\begin{align*}
 I_{el} &\geq  \frac{e}{\permitELscalar\boltz\temp} \sum_l\Sign(\zl)\scp{\chargeELold}{(\abs{\zl}\cl)^2}_\Omega
               -\delta \sum_l\abs{\zl}\norm{\grad\cl}{\Lp{2}{}}^2 \\
        &~~~~ -\frac{2e^2\max_l\abs{\zl}^2}{(\permitELscalar\boltz\temp)^2}\sqbrac{\delta^{-1}\norm{\sigma}{\Lp[\GammaT]{\infty}{}}^2+\norm{\rho_b}{\Lp[\OmegaT]{\infty}{}}} \sum_l\abs{\zl}\norm{\cl}{\Lp{2}{}}^2 \\
        &~~~~ := \labelInt + I.b + I.c.
\end{align*}
Now it remains to control the integral~$\labelInt$. As we assumed $\cl=\clold$, the free charge density~$\chargeELold$ is given by $\chargeELold=\chargeELlongTwo$.
This shows with \cref{lemma:sign}
\label{page:signCondition}
\begin{align*}
 \labelInt = \frac{e}{\permitELscalar\boltz\temp}~\scp{\chargeELlongTwo}{(\zl[1]\cl[1])^2-(\abs{\zl[2]}\cl[2])^2}_\Omega \geq 0.
\end{align*}
This sign condition%
\footnote{Thanks to $\labelInt\geq0$, we can avoid to bound the integral~$\labelInt$ by suitable norms of its integrands, which would cause serious problems. See \cref{sec:Introduction} for further details.}
only holds true, since we use the weighted test function~$\test=\abs{\zl}\cl$ instead of $\test=\cl$.
% Substituting $\labelInt\geq0$ into the above intermediate estimate, leads for electric drift integral to
% \begin{align*}
%  I_{el}&\geq -\delta \sum_l\abs{\zl}\norm{\grad\cl}{\Lp{2}{}}^2
%              -\frac{2e^2}{(\permitELscalar\boltz\temp)^2}\sqbrac{\delta^{-1}\norm{\sigma}{\Lp[\GammaT]{\infty}{}}^2+\norm{\rho_b}{\Lp[\OmegaT]{\infty}{}}} \sum_l\abs{\zl}\norm{\cl}{\Lp{2}{}}^2 .
% \end{align*}
For the surface integrals, we involve \ref{Assump:BoundData}, \cref{lemma:interpol-boundary}, and \Young. Thereby, we get
\begin{align*}
 &      \sum_l \scp{\gl}{\abs{\zl}\cl}_\Omega
 ~\leq~ \sum_l\abs{\zl}\norm{\cl}{\Lp[\Gamma]{2}{}}^2 + \max_l\abs{\zl}\sum_l\norm{\gl}{\Lp[\Gamma]{2}{}}^2 \\
 &\leq  \delta \sum_l \abs{\zl} \norm{\grad\cl}{\Lp{2}{}}^2
        +2\delta^{-1} \sum_l \abs{\zl}\norm{\cl}{\Lp{2}{}}^2
        +\max_l\abs{\zl} \sum_l\norm{\gl}{\Lp[\Gamma]{2}{}}^2 ~.
\end{align*}
Applying \ref{Assump:Reaction}, \Young, and recalling $\cl\geq0$, results for the reaction integrals in
\begin{align*}
  &    \sum_l\scp{\theta\Rl(\Cl)}{\abs{\zl}\cl}_\Omega
   \leq \theta \max_l C_{\Rl} \scp{\abs{\Cl}}{\sum_l\abs{\zl}\cl}_\Omega
   \leq \theta \max_l C_{\Rl} \scp{\sum_l\abs{\zl}\cl}{\sum_l\abs{\zl}\cl}_\Omega \\
  &\leq 3\theta \max_l\abs{\zl}C_{\Rl} \sum_l \norm{\abs{\zl}\cl}{\Lp{2}{}}^2 ~.
\end{align*}
By combining the preceding estimates, we deduce with the choice~$\delta:=\frac{\alpha_D}{6}$ the estimate
\begin{align}\label{eq:energy-energy1}
 &    \frac{\theta}{2}\derr\sum_l \abs{\zl}\norm{\cl}{\Lp{2}{}}^2 + \frac{\alpha_D}{2} \sum_l\abs{\zl}\norm{\grad\cl}{\Lp{2}{}}^2 \nonumber\\
 &\leq  \frac{2e^2}{(\permitELscalar\boltz\temp)^2}\sqbrac{\frac{6}{\alpha_D} \norm{\sigma}{\Lp[\GammaT]{\infty}{}}^2+\norm{\rho_b}{\Lp[\OmegaT]{\infty}{}}} \sum_l\abs{\zl}\norm{\cl}{\Lp{2}{}}^2  \nonumber\\
 &~~~~  +\sqbrac{\frac{6}{\alpha_D} \norm{f}{\Lp[\GammaT]{\infty}{}}^2 +  3\theta \max_l\abs{\zl}C_{\Rl} }\sum_l\abs{\zl}\norm{\cl}{\Lp{2}{}}^2
        +\max_l\abs{\zl}\sum_l\norm{\gl}{\Lp[\Gamma]{2}{}}^2~.
\end{align}
For ease of readability, we introduce the abbreviation
\begin{align*}
 B_0 &:= \max\!\brac{\frac{2}{\theta},\frac{2}{\alpha_D}} \frac{12e^2\max_l\abs{\zl}^2}{\alpha_D(\permitELscalar\boltz\temp)^2}
         \left[ \norm{\sigma}{\Lp[\GammaT]{\infty}{}}^2 + \norm{\rho_b}{\Lp[\OmegaT]{\infty}{}}  + \norm{f}{\Lp[\GammaT]{\infty}{}}^2 +  3\theta \max_lC_{\Rl}  \right].
\end{align*}
Thus, we immediately obtain from \eqref{eq:energy-energy1} together with \Grownwall\
\begin{align*}
      \sum _l\abs{\zl}\norm{\cl}{\Lp[I]{\infty}{;\Lp{2}{}}}^2
 \leq  e^{B_0\timeEnd} \sqbrac{\sum_l\abs{\zl}\norm{\clstart}{\Lp{2}{}}^2 +  \max_l\abs{\zl}\sum_l\norm{\gl}{\Lp{2}{}}^2}
 ~:=   \hat{C}^2_0(\timeEnd).
\end{align*}
We substitute this bound into \eqref{eq:energy-energy1} and we integrate in time over $[0,\timeEnd]$. This yields
\begin{align}\label{eq:energy-statedbound-1}
  &     \sum _l \abs{\zl}\sqbrac{\norm{\cl}{\Lp[I]{\infty}{;\Lp{2}{}}} + \norm{\grad\cl}{\Lp[\OmegaT]{2}{}} } \nonumber\\
  &\leq \hat{C}_0 + B_0^{\frac{1}{2}}\max_l\abs{\zl}\timeEnd^{\frac{1}{2}}\hat{C}_0+  \max_l\abs{\zl}\sum_l\norm{\gl}{\Lp{2}{}}
   ~:=  C_0(\timeEnd)~.
\end{align}
%
%
%=======================================================
% ENERGY - CASE 2: \cl\neq\clold
%=======================================================
\underline{Case 2: $\cl\neq\clold$}~~
Again, we test equations~\eqref{eq:transportWeak-fpi} with $\test:=\abs{\zl}\cl\in\spaceTC$. The above estimates for the respective integrals remain unchanged,
except for the integral~$\labelInt$, which does not fulfill a sign condition this time. Thus, we now bound this integral with \Holder\ by
\begin{align*}
 \labelInt &=    \frac{e}{\permitELscalar\boltz\temp} \scp{\chargeELold}{\sum_l\Sign(\zl)\brac{\abs{\zl}\cl}^2}_\Omega
           ~\leq \frac{e}{\permitELscalar\boltz\temp} \max_l\abs{\zl}^2\norm{\Clold}{\Lp{\infty}{}}~~ \sum_l\abs{\zl}\norm{\cl}{\Lp{2}{}}^2.
\end{align*}
Herein, the constant depends on the $L^\infty$-norms of the $\clold$. However, this is uncritical in this case as we assumed $\cl\neq\clold$. Furthermore, in \cref{def:FixOp}
we introduced the space~$X$ and the set~$\fixSet\subset X$. Furthermore, we supposed $\Clold\in\fixSet$, which ensures that the $L^\infty$-norms of the $\clold$ remain finite.
Thus, provided we know $\norm{\Clold}{\Lp{\infty}{}}\leq R$ for all $\Clold\in\fixSet$, the constant in the above estimate just depends on an additional parameter~$R$.
In conclusion, with the redefined constant
\begin{align*}
 B_0 &:=\frac{24e^2\max_l\abs{\zl}^2C_{\Rl}}{\min(\theta,\alpha_D)\alpha_D(\permitELscalar\boltz\temp)^2}
         \left[ \norm{\sigma}{\Lp[\GammaT]{\infty}{}}^2 + \norm{\rho_b}{\Lp[\OmegaT]{\infty}{}}  + \norm{f}{\Lp[\GammaT]{\infty}{}}^2 +  3\theta + \norm{\Clold}{\Lp{\infty}{}} \right],
\end{align*}
we obtain analogously to \eqref{eq:energy-statedbound-1}
\begin{align}\label{eq:energy-statedbound-2}
       \sum _l \abs{\zl}\sqbrac{\norm{\cl}{\Lp[I]{\infty}{;\Lp{2}{}}} + \norm{\grad\cl}{\Lp[\OmegaT]{2}{}} }
 \leq  C_0(\timeEnd,\norm{\Clold}{\Lp{\infty}{}} )~.
\end{align}
\end{proof}
\begin{remark}\label{remark:boundForFreeCharge}
 In particular, \cref{lemma:energy} ensures $\norm{\chargeELold}{\Lp[I]{\infty}{;\Lp{2}{}}} \leq \max_l\abs{\zl} C_0$.
 This uniform bound holds for both cases, $\chargeEL=\chargeELold$ and $\chargeEL\neq\chargeELold$.
 \hfill$\square$
\end{remark}
Next, we show that the chemical species~$\cl$ are bounded. As the proof is rather long and technical, we separate this proof from the proof of the remaining a~priori bounds in \cref{thm:aprioriBounds}.
%
%
%============================
%============================
% LEMMA: Moser Bound
%============================
%============================
\begin{lemma}[Boundedness]\label{lemma:bounded}
Let \ref{Assump:Geom}--\ref{Assump:back-charge} be valid and let $\brac{\fieldEL,\potEL,\fieldF,\press,\Cl}\in \setR^{4+2n}$ be a weak solution
of \eqref{eq:Model-1a}--\eqref{eq:Model-1j} according to \cref{def:FixOp}. Then, we have
\begin{align*}
  \sum_l \norm{\cl}{\Lp[\OmegaT]{\infty}{}} ~\leq~ C_M~.
\end{align*}
Herein, the dependency of the constant is
\begin{align*}
 C_M=C_M \!\brac{\timeEnd,\max_l\abs{\zl},\norm{\gl}{\Lp[\OmegaT]{\infty}{}},\norm{f}{\Lp[\GammaT]{\infty}{}},\norm{\sigma}{\Lp[\GammaT]{\infty}{}},\norm{\rho_b}{\Lp[\OmegaT]{\infty}{}},\norm{\clstart}{\Lp{\infty}{}} } .
\end{align*}
\end{lemma}
\begin{proof}
As we have already established a lower bound for $\cl$ in \cref{lemma:nonnegative}, it remains to show an upper bound.
To this end, we subsequently apply Moser's iteration technique, cf. \cite{Moser1,Moser2}. More precisely, we follow the
proof of \cite[Theorem 6.15]{Liebermann-book} with a modified test function.
\par
Henceforth, we use the truncated solutions  $\clm:=\min(\cl,m)$. For ease of readability, we split the proof into several steps.\\[2.0mm]
%
%
%=======================================================
% MOSER: STEP 1 - preliminary energy estimate
%=======================================================
\underline{Step 1: preliminary energy estimates}\\
The crucial step in Moser's iteration technique is to derive an energy estimate for $\cl$ to arbitrary high powers, i.e., for $\cl^{\alpha+1}$ with $\alpha\geq0$.
For that purpose, we test equations~\eqref{eq:transportWeak-fpi} by $\varphi := \brac{\clm}^{2\alpha+1}$ for $\alpha\geq0$,
we sum over $l=1,2$, and we bound the respective integrals. This part of the proof is related to the proof of \cref{lemma:energy}. For this reason, we just briefly repeat the similar parts.
Firstly, using the above test function yields with \ref{Assump:Ellip} for the diffusion integrals
\begin{align*}
      (2\alpha+1) \sum_l\scp{\Dl\grad\cl}{(\clm)^{2\alpha}\grad\clm}_\Omega
 \geq \frac{\alpha_D}{\alpha+1} \sum_l\norm{\grad(\clm)^{\alpha+1}}{\Lp{2}{}}^2~.
\end{align*}
The convection integrals, we firstly transform with integration by parts and inserting equation \eqref{eq:darcyWeakDIV-fpi}. Then, we involving \Holder\ and \cref{lemma:interpol-boundary}
(with a rescaled parameter $\delta$). Thereby, we arrive for the convection integrals at to
\begin{align*}
 &     -\sum_l\scp{{\cl}\fieldF}{\grad\brac{\clm}^{2\alpha+1}}_\Omega
  =    -\sum_l\frac{2\alpha+1}{(2\alpha+2)} \scp{f}{\brac{\clm}^{2\alpha+2}}_\Gamma \\
 &\geq -\delta \sum_l\norm{\grad\brac{\clm}^{\alpha+1}}{\Lp{2}{}}^2
  ~-~  2\delta^{-1} \norm{f}{\Lp[\GammaT]{\infty}{}}^2 \sum_l\norm{\brac{\clm}^{\alpha+1}}{\Lp{2}{}}^2.
\end{align*}
Analogously, we transform the electric drift integrals with equation~\eqref{eq:gaussWeakDIV-fpi} to
\begin{align*}
 I_{el}:=&  -\frac{e}{\permitELscalar\boltz\temp} \sum_l\scp{\zl{\cl}\fieldEL}{\grad\brac{\clm}^{2\alpha+1}}_\Omega
         ~= -\frac{e(2\alpha+1)}{\permitELscalar\boltz\temp(2\alpha+2)} \sum_l\zl\scp{\fieldEL}{\grad\brac{\clm}^{2\alpha+2}}_\Omega\\
         &= \frac{e(2\alpha+1)}{\permitELscalar\boltz\temp(2\alpha+2)} \sum_l\zl\sqbrac{~\scp{\rho_b+\chargeELold}{\brac{\clm}^{2\alpha+2}}_\Omega - \scp{\sigma}{\brac{\clm}^{2\alpha+2}}_\Gamma~}.
\end{align*}
Next, we apply \Holder, and \cref{lemma:interpol-boundary}. Thereby, we come (with a rescaled $\delta$) to
\begin{align*}
 I_{el} &\geq  \frac{e(2\alpha+1)\max_l\abs{\zl} }{\permitELscalar\boltz\temp(2\alpha+2)}\norm{\chargeELold}{\Lp{2}{}}\sum_l\norm{\brac{\clm}^{2\alpha+2}}{\Lp{2}{}}
               -\delta \sum_l\norm{\grad\brac{\clm}^{\alpha+1}}{\Lp{2}{}}^2 \\
        &~~~~  -\underbrace{\frac{2e^2}{(\permitELscalar\boltz\temp)^2}\max_l\abs{\zl}\sqbrac{\delta^{-1}\norm{\sigma}{\Lp[\GammaT]{\infty}{}}^2+\norm{\rho_b}{\Lp[\OmegaT]{\infty}{}}} }_{=:K_0} \sum_l\norm{\brac{\clm}^{\alpha+1}}{\Lp{2}{}}^2 \\
        &~~~~  := \labelInt + I.b + I.c.
\end{align*}
We are done for integrals~$I.b$ and $I.c$. As to the integral~$\labelInt$, we apply \GagNirenberg, cf. \cite{nirenberg1959}, \Young, \cref{lemma:energy}, and \cref{remark:boundForFreeCharge}. This shows
\label{page:withoutSignCondition}
\begin{align*}
 \labelInt &=     \frac{e(2\alpha+1)}{\permitELscalar\boltz\temp(2\alpha+2)}\max_l\abs{\zl} \norm{\chargeELold}{\Lp{2}{}}\sum_l\norm{\brac{\clm}^{\alpha+1}}{\Lp{4}{}}^2  \\
           &\geq -\frac{e}{\permitELscalar\boltz\temp} \max_l\abs{\zl}^2 C_0C_{gn} \sum_l\norm{\brac{\clm}^{\alpha+1}}{\Lp{2}{}}^{(4-n)/2}\norm{\brac{\clm}^{\alpha+1}}{\Hk{1}{}}^{n/2}\\
           &\geq -\delta \sum_l\norm{\grad\brac{\clm}^{\alpha+1}}{\Lp{2}{}}^2 - \underbrace{C(\delta^{-1})\frac{e^4}{(\permitELscalar\boltz\temp)^4} \max_l\abs{\zl}^8 C_0^4C_{gn}^4}_{=:K_1} \sum_l\norm{\brac{\clm}^{\alpha+1}}{\Lp{2}{}}^2 .
\end{align*}
Substituting this into the above intermediate estimate, results for electric drift integral in
\begin{align*}
  I_{el} \geq -\delta \sum_l\norm{\grad\brac{\clm}^{\alpha+1}}{\Lp{2}{}}^2
              -\sqbrac{K_0+K_1} \sum_l\norm{\brac{\clm}^{\alpha+1}}{\Lp{2}{}}^2 .
\end{align*}
For the surface integrals, we involve \ref{Assump:BoundData}, \cref{lemma:interpol-boundary}. Furthermore, we apply \Young\ with $q=\frac{2\alpha+2}{2\alpha+1}$, $p=2\alpha+2$. Thereby, we get
\begin{align*}
 &      \sum_l \scp{\gl}{\brac{\clm}^{2\alpha+1}}_\Omega
 ~\leq~ \sum_l\norm{\brac{\clm}^{2\alpha+2}}{\Lp[\Gamma]{1}{}} + \sum_l\norm{\gl^{2\alpha+2}}{\Lp[\Gamma]{1}{}} \\
 &\leq  \delta \sum_l \norm{\grad\brac{\clm}^{\alpha+1}}{\Lp{2}{}}^2
        +2\delta^{-1} \sum_l \norm{\brac{\clm}^{\alpha+1}}{\Lp{2}{}}^2
        +\sum_l\norm{\gl^{\alpha+1}}{\Lp[\Gamma]{2}{}}^2 ~.
\end{align*}
A combination of the preceding estimates, together with the not yet considered time integrals and reaction integrals, yields with the choice~$\delta:=\frac{\alpha_D}{6(\alpha+1)}$
the preliminary energy estimate
\begin{align}\label{eq:moser-energy-preliminary-1}
 &     \theta \sum _l\dualp{\dert\cl}{\brac{\clm}^{2\alpha+1}}_{1,\Omega} + \frac{\alpha_D}{2(\alpha+1)}~\sum_l\norm{\grad\brac{\clm}^{\alpha+1}}{\Lp{2}{}}^2 \nonumber\\
 &\leq  \frac{12(\alpha+1)}{\alpha_D}\sqbrac{\norm{f}{\Lp[\GammaT]{\infty}{}}+1+K_0+K_1}\sum_l\norm{\brac{\clm}^{\alpha+1}}{\Lp{2}{}}^2 \nonumber\\
 &~~~~  +\sum_l\norm{\gl^{\alpha+1}}{\Lp[\Gamma]{2}{}}^2 + \sum_l\scp{\theta\Rl(\Cl)}{\brac{\clm}^{2\alpha+1}}_\Omega
\end{align}
%
%
%===========================================
% MOSER: STEP 2 - base case
%===========================================
\underline{Step 2: base case:}~~
Before we continue the proof, we define for $j\in\setN_0$ a sequence of exponents~$\alpha_j$ by
\begin{align} \label{eq:moser-exponents}
1+\alpha_j &:= \brac{\frac{n+2}{n}}^j ~~~\Hence~ 1+\alpha_j = (1+\alpha_{j-1})(1+\alpha_1)~.
\end{align}
Additionally, we cite from \cite[Proposition~3.2]{DiBenedetto-book} the parabolic embedding
\begin{align}\label{eq:moser-parab-embedding}
 \norm{\cl}{\Lp[\OmegaT]{2\frac{n+2}{n}}{}} \leq C_S\!\brac{\timeEnd,\Gamma,n}~\sqbrac{ \norm{\cl}{\Lp[I]{\infty}{;\Lp{2}{}}} + \norm{\grad\cl}{\Lp[\OmegaT]{2}{}} ~}.
\end{align}
We now begin Moser's iteration procedure, which we rigorously formulate as mathematical induction.
% The idea behind Moser's iteration procedure is as follows:
% Provided that $\cl^{\alpha_j}\in\Lp[\OmegaT]{2}{}$, we firstly establish an energy estimate for $\cl^{\alpha_j}$. Secondly, we conclude with \eqref{eq:moser-parab-embedding}
% that even $\cl^{\alpha_{j+1}}$ lies in $\Lp[\OmegaT]{2}{}$. Since $\cl^{\alpha_j}\in\Lp[\OmegaT]{2}{}$ is equivalent to $\cl\in\Lp[\OmegaT]{2\alpha_j}{}$, we successfully increase the
% integrability exponent of $\cl$ from $2\alpha_j$ to $2\alpha_{j+1}$.
% This iteration procedure leads with $\alpha_j \rightarrow\infty$ as $j\rightarrow\infty$ finally to $\cl\in\Lp[\OmegaT]{p}{}$ for all $p\in[2,\infty]$.
% \par
% Let us start the mathematical induction with $j=0$.
We start with $j=0$, Thus, we have  $1+\alpha_0=1$, which means $\alpha_0=0$. Substituting this into \eqref{eq:moser-energy-preliminary-1},
we can safely let $m\rightarrow\infty$, as every integral remains finite. Furthermore, the reaction integrals and the time integrals, we estimate exactly as
in the proof of \cref{lemma:energy}. Thereby, we rediscover a slightly modified version of \eqref{eq:energy-energy1}. I.e., we arrive with $m\rightarrow\infty$ at
\begin{align*}
 &     \frac{\theta}{2} \sum _l\derr\norm{\cl}{\Lp{2}{}}^2  + \frac{\alpha_D}{2} \sum_l\norm{\grad\cl}{\Lp{2}{}}^2 \nonumber\\
 &\leq  \frac{12}{\alpha_D} \sqbrac{ \norm{f}{\Lp[\GammaT]{\infty}{}}+1+K_0+K_1+3\theta\max_l C_{\Rl} }\sum_l\norm{\cl}{\Lp{2}{}}^2 + \sum_l\norm{\gl}{\Lp[\Gamma]{2}{}}^2~.
\end{align*}
Again, we deduce with \Grownwall\
\begin{subequations}
\begin{align}\label{eq:moser-bound1-0}
  &     \sum _l \sqbrac{\norm{\brac{\cl}}{\Lp[I]{\infty}{;\Lp{2}{}}} + \norm{\grad\brac{\cl}}{\Lp[\OmegaT]{2}{}} }
   \leq \hat{C}_0 + B_0^{\frac{1}{2}}\timeEnd^{\frac{1}{2}}\hat{C}_0+ \sum_l\norm{\gl}{\Lp{2}{}}
   ~:=  C_0(\timeEnd)~.
\end{align}
This time, we have denoted the constants by
\begin{align*}
          B_0:=&\min\!\brac{\frac{2}{\theta}, \frac{2}{\alpha_D}} \frac{12\max_l\abs{\zl}}{\alpha_D}\sqbrac{\norm{f}{\Lp[\GammaT]{\infty}{}}+1 + K_0+K_1 + 3\theta \max_l C_{\Rl} },\\
 \hat{C}^2_0 := & e^{B_0\timeEnd} \sqbrac{\sum_l\norm{\clstart}{\Lp{2}{}}^2 + \sum_l\norm{\gl}{\Lp{2}{}}^2} .
\end{align*}
By combining \eqref{eq:moser-bound1-0}, the definition of $\alpha_1$ in \eqref{eq:moser-exponents}, and the embedding~\eqref{eq:moser-parab-embedding}, we finally obtain
(with the elementary inequality $(a+b)^{1/p} \leq a^{1/p}+b^{1/p}$ for $a,b\geq0$ and $p\geq 1$ )
\begin{align}\label{eq:moser-bound2-0}
 &      \sqbrac{\sum _l \norm{\cl}{\Lp[\OmegaT]{2(1+\alpha_1)}{}}^{2(1+\alpha_1)}}^{\frac{1}{2(1+\alpha_1)}}
  \leq  \sum _l \norm{\cl}{\Lp[\OmegaT]{2(1+\alpha_1)}{}}
  =     \sum _l \norm{\cl}{\Lp[\OmegaT]{2\frac{n+2}{n}}{}}  \nonumber\\
 &\leq C_S\sum _l \sqbrac{\norm{\cl}{\Lp[I]{\infty}{;\Lp{2}{}}} + \norm{\grad\cl}{\Lp[\OmegaT]{2}{}} }
  \leq C_SC_0~.
\end{align}
\end{subequations}
%
%
%===========================================
% MOSER: STEP 3 - inductive hypothesis
%===========================================
\underline{Step 3: induction hypothesis}~~
Let $j\in\setN_0$ and let $\alpha_j$ be the corresponding exponent defined in \eqref{eq:moser-exponents}. We suppose, there exists a constant $C_{j-1}(\timeEnd)$ such that we have
\begin{subequations}
\begin{align}
 \sum _l \sqbrac{\norm{\brac{\cl}^{1+\alpha_{j-1}}}{\Lp[I]{\infty}{;\Lp{2}{}}} + \norm{\grad\brac{\cl}^{1+\alpha_{j-1}}}{\Lp[\OmegaT]{2}{}} } &\leq C^{1+\alpha_{j-1}}_{j-1}(\timeEnd), \label{eq:moser-bound1-j} \\
      \sum _l \norm{\brac{\cl}^{\alpha_j}}{\Lp[\OmegaT]{2}{}} &\leq C_{j-1}(\timeEnd)~. \label{eq:moser-bound2-j}
\end{align}
\end{subequations}
We note that \eqref{eq:moser-bound1-j} and \eqref{eq:moser-bound2-j} reduce for $j=0$ exactly to \eqref{eq:moser-bound1-0} and \eqref{eq:moser-bound2-0}.\\[2.0mm]
%
%===========================================
% MOSER: STEP 4 - inductive step
%===========================================
\underline{Step 4: inductive step:}~~
We return with the choice of $\alpha:=\alpha_j$ to \eqref{eq:moser-energy-preliminary-1} and we can safely let $m\rightarrow\infty$, as we already know \eqref{eq:moser-bound2-j}.
For the reaction integrals, this leads analogously to the base case with \Young\ $(p=\frac{2\alpha_j+1}{2\alpha_j+2}, q=2\alpha_j+2)$ to
\begin{align*}
  \sum_l\scp{\theta\Rl(\Cl)}{\brac{\cl}^{2\alpha_j+1}}_\Omega \leq 3\theta \max_l C_{\Rl} \sum_l \norm{\brac{\cl}^{\alpha_j+1}}{\Lp{2}{}}^2 ~.
\end{align*}
Additionally, we estimate the time integrals analogously to the base case. Thus, after letting $m\rightarrow\infty$, by incorporating the bounds for the reaction integrals and the time integrals,
and by introducing the abbreviations
\begin{align*}
 A_j:=&\min\!\brac{ \frac{\theta}{(2\alpha_j+2)}, \frac{\alpha_D}{2(\alpha_j+1)} },\\
 B_j:=&\frac{12(\alpha_j+1)}{A_j\alpha_D}\sqbrac{\norm{f}{\Lp[\GammaT]{\infty}{}}+ 1 + K_0 + K_1 3\theta \max_l C_{\Rl} },
\end{align*}
we finally obtain the energy estimate
\begin{align*}
  \derr\!\sum_l\norm{\brac{\cl}^{\alpha_j+1}}{\Lp{2}{}}^2  \!+\! \sum_l\norm{\grad\brac{\cl}^{\alpha_j+1}}{\Lp{2}{}}^2
   \!\leq\! B_j\! \sum_l\norm{\brac{\cl}^{\alpha_j+1}}{\Lp{2}{}}^2 \!\!+\! \sum_l\norm{\gl^{\alpha_j+1}}{\Lp[\Gamma]{2}{}}^2.
\end{align*}
Hence, we conclude with \Grownwall\ the uniform bound
\begin{align*}
      \sum _l\norm{\brac{\cl}^{\alpha_j+1}}{\Lp[I]{\infty}{;\Lp{2}{}}}^2
 \leq \underbrace{e^{B_j\timeEnd} \sqbrac{\sum_l\norm{\brac{\clstart}^{\alpha_j+1}}{\Lp{2}{}}^2 + \sum_l\norm{\gl^{\alpha_j+1}}{\Lp{2}{}}^2} }_{:=  \hat{C}^{2(\alpha_j+1)}_j(\timeEnd)}.
\end{align*}
Next, we integrate the above energy estimate in time over $[0,\timeEnd]$ and we involve the preceding bound. Thereby, we deduce the stated inequality~\eqref{eq:moser-bound1-j}
\begin{subequations}
\begin{align}\label{eq:moser-bound1-j+1}
  &     \sum _l \sqbrac{\norm{\brac{\cl}^{\alpha_j+1}}{\Lp[I]{\infty}{;\Lp{2}{}}} + \norm{\grad\brac{\cl}^{\alpha_j+1}}{\Lp[\OmegaT]{2}{}} } \nonumber\\
  &\leq \hat{C}^{\alpha_j+1}_j + B_j^{\frac{1}{2}}\timeEnd^{\frac{1}{2}}\hat{C}^{\alpha_j+1}_j+ \sum_l\norm{\gl^{\alpha_j+1}}{\Lp{2}{}}
   ~:=  C^{\alpha_j+1}_j(\timeEnd)~.
\end{align}
Furthermore, with the definition of $\alpha_j$ in \eqref{eq:moser-exponents}, and the embedding~\eqref{eq:moser-parab-embedding}, we arrive at the stated bound~\eqref{eq:moser-bound2-j}
\begin{align}\label{eq:moser-bound2-j+1}
 &      \sqbrac{ \sum _l \norm{\cl}{\Lp[\OmegaT]{2(1+\alpha_{j+1})}{}}^{2(1+\alpha_{j+1})} }^{\frac{1}{2(1+\alpha_{j+1})}}
  \leq  \sqbrac{ \sum _l \norm{\brac{\cl}^{\alpha_j+1}}{\Lp[\OmegaT]{2(1+\alpha_1)}{}} }^{\frac{1}{(1+\alpha_j)}} \nonumber\\
 &\leq  \sqbrac{ C_S\sum _l \brac{\norm{\brac{\cl}^{\alpha_j+1}}{\Lp[I]{\infty}{;\Lp{2}{}}} + \norm{\grad\brac{\cl}^{\alpha_j+1}}{\Lp[\OmegaT]{2}{}} } }^{\frac{1}{(1+\alpha_j)}} \nonumber\\
 &\leq  \sqbrac{C_SC^{1+\alpha_j}_j}^{\frac{1}{(1+\alpha_j)}}
 ~\leq  C_S^{\frac{1}{(1+\alpha_j)}}C_j ~.
\end{align}
This shows that the induction hypothesis holds for all $j\in\setN$.\\[2.0mm]
%
%
%
%===========================================
% MOSER: STEP 5 - limit case
%===========================================
\underline{Step 5: limit case:}~~
We now consider the limit case $j\rightarrow\infty$. First of all, we note that we have
 \begin{align*}
  \norm{\Cl}{\Lp[\OmegaT]{p}{}}      &:= \brac{ \sum_l \norm{\brac{\cl}}{\Lp[\OmegaT]{p}{}}^p }^{\frac{1}{p}} && 1\leq p < \infty~, \\
  \norm{\Cl}{\Lp[\OmegaT]{\infty}{}} &:= \sum_l \norm{\brac{\cl}}{\Lp[\OmegaT]{\infty}{}}~,                   && p=\infty~.
 \end{align*}
Hence, we can rewrite \eqref{eq:moser-bound2-j+1} as
\begin{align}\label{eq:moser-bound2-j+1-new}
       \norm{\Cl}{\Lp[\OmegaT]{2(1+\alpha_{j+1})}{}}
 =     \sqbrac{ \sum _l \norm{\cl}{\Lp[\OmegaT]{2(1+\alpha_{j+1})}{}}^{2(1+\alpha_{j+1})} }^{\frac{1}{2(1+\alpha_{j+1})}}
 \leq  C_S^{\frac{1}{(1+\alpha_j)}}C_j ~.
\end{align}
\end{subequations}
Our goal is to show that this inequality holds even in the limit case $j=\infty$. For that purpose, we recall that $a^{1/p}\rightarrow1$ as $p\rightarrow\infty$ for all $a>0$.
Furthermore, from the definition of $\alpha_j$ in \eqref{eq:moser-exponents}, we know that $\alpha_j\rightarrow\infty$ as $j\rightarrow\infty$.
Thus, with the definition of the constants $C_j$, we obtain in the limit
\begin{align*}
C_\infty &:=   \lim_{j\rightarrow\infty} \brac{C_S^{\frac{1}{(1+\alpha_j)}}C_j}
         ~=    \brac{\lim_{j\rightarrow\infty} C_S^{\frac{1}{(1+\alpha_j)}} } \brac{\lim_{j\rightarrow\infty} C_j }
         ~=    \lim_{j\rightarrow\infty} C_j \\
         &\leq \lim_{j\rightarrow\infty} \hat{C}_j + \lim_{j\rightarrow\infty} B_j^{\frac{1}{2(1+\alpha_j)}} \timeEnd^{\frac{1}{2(1+\alpha_j)}} \hat{C}_j + \lim_{j\rightarrow\infty}\sum_l\norm{\gl^{\alpha_j+1}}{\Lp{2}{}}^{\frac{1}{(1+\alpha_j)}}\\
         &:=   C_{\infty,1} + C_{\infty,2} + C_{\infty,3}.
\end{align*}
Next, we recall from \cite[Theorem~2.14]{Adams2-book} that for all $u\in\Lp{\infty}{}$ holds
\begin{align}\label{eq:LpNorm->MaxNorm}
 \norm{u}{\Lp{\infty}{}} = \lim\limits_{p\rightarrow\infty} \norm{u}{\Lp{p}{}}~.
\end{align}
This reveals immediately with \ref{Assump:BoundData} that $C_{\infty,3}=\lim_{j\rightarrow\infty} \sum_l\norm{\gl}{\Lp{2(\alpha_j+1)}{}} = \sum_l\norm{\gl}{\Lp{\infty}{}}$.
Furthermore, we arrive with \eqref{eq:LpNorm->MaxNorm}, the definition of $\hat{C}_j$ and $(e^x)^y=e^{xy}$ for $x,y>0$ at
\begin{align*}
 C_{\infty,1} &\leq e^{\frac{B_0\timeEnd}{2}} \brac{\sum_l\norm{\clstart}{\Lp{\infty}{}} +  \sum_l\norm{\gl}{\Lp{\infty}{}}},
\end{align*}
where we used $\lim_{j\rightarrow\infty} \exp\!\brac{\frac{B_j\timeEnd}{2(1+\alpha_j)} } \leq \exp\!\brac{\frac{B_0\timeEnd}{2} }$ due to the definition of $B_j$ and $B_0$.
Analogously, we get for $C_{\infty,2}$ with $j^{1/j}\rightarrow1$ as $j\rightarrow\infty$
\begin{align*}
 C_{\infty,2} &=   \lim_{j\rightarrow\infty} B_j^{\frac{1}{2(1+\alpha_j)}} \timeEnd^{\frac{1}{2(1+\alpha_j)}} \hat{C}_j
              \leq e^{\frac{B_0\timeEnd}{2}} \brac{\sum_l\norm{\clstart}{\Lp{\infty}{}} +  \sum_l\norm{\gl}{\Lp{\infty}{}}}.
\end{align*}
Combining the preceding estimates, shows that
\begin{align*}
C_\infty \leq 2 e^{\frac{B_0\timeEnd}{2}} \brac{\sum_l\norm{\clstart}{\Lp{\infty}{}} +  \sum_l\norm{\gl}{\Lp{\infty}{}}} + \sum_l\norm{\gl}{\Lp{\infty}{}}.
\end{align*}
We now can safely let $j\rightarrow\infty$ in \eqref{eq:moser-bound2-j+1-new}. Thereby we finally arrive together with \eqref{eq:LpNorm->MaxNorm} at
\begin{align*}
 &\sum_l \norm{\cl}{\Lp[\OmegaT]{\infty}{}} = \norm{\Cl}{\Lp[\OmegaT]{\infty}{}} = \lim_{j\rightarrow\infty} \norm{\Cl}{\Lp[\OmegaT]{2(\alpha_j+1)}{}} \\
 &\leq 2 e^{\frac{B_0\timeEnd}{2}} \brac{\sum_l\norm{\clstart}{\Lp{\infty}{}} +  \sum_l\norm{\gl}{\Lp{\infty}{}}} + \sum_l\norm{\gl}{\Lp{\infty}{}}
 ~=:   C_M(\timeEnd)~.
\end{align*}
\end{proof}
Finally, we show the desired a~priori bounds for a solution of the \dpnp. These a~priori bounds are crucial for the proof of \cref{thm:existence}.
%
%
%============================
%============================
% THEOREM: A priori Bounds
%============================
%============================
\begin{thm}[A priori Bounds]\label{thm:aprioriBounds}
Let \ref{Assump:Geom}--\ref{Assump:back-charge} be valid and let $\brac{\fieldEL,\potEL,\fieldF,\press,\Cl}\in \setR^{4+2n}$ be a weak solution
of \eqref{eq:Model-1a}--\eqref{eq:Model-1j} according to \cref{def:weaksolution}. Then, we have
\begin{align*}
 & \norm{\potEL}{\Lp[I]{\infty}{;\Lp{2}{}}} + \norm{\fieldEL}{\Lp[I]{\infty}{;\Lp{2}{}}} ~\leq~ C(\timeEnd)~, \\
 & \norm{\press}{\Lp[I]{\infty}{;\Lp{2}{}}} + \norm{\fieldF}{\Lp[I]{\infty}{;\Lp{2}{}}}  ~\leq~ C(\timeEnd)~, \\
 & \sum_l\sqbrac{ \norm{\cl}{\Lp[I]{\infty}{;\Lp{2}{}}}+\norm{\cl}{\Lp[I]{2}{;\spaceTC}}+\norm{\cl}{\Hk[I]{1}{;\spaceTC^\ast}}+\norm{\cl}{\Lp[\OmegaT]{\infty}{}} }~\leq~ C(\timeEnd)~.
\end{align*}
\end{thm}
\begin{remark}
\textnormal{
 In the following proof, we derive the constants of the stated a~priori estimates in detail. This reveals how the constants depend on the data.
However, for the following \cref{thm:existence}, especially the dependency of the end time~$\timeEnd$ is of interest. This is the reason why we just stated $C=C(\timeEnd)$ for the constants .
}\hfill$\square$
\end{remark}
\begin{proof} For ease of readability, we split the proof of the stated a~priori estimates into several steps.\\[2.0mm]
%
%====================================================================
% Step 1.1: NERNST-PLANCK - energy-bound for cl=\clold
%====================================================================
\underline{Step~1.1 -- energy estimates and boundedness for $\cl$:}~~
These a~priori bounds are shown in \cref{lemma:energy} and \cref{lemma:bounded}.\\[2.0mm]
%
%====================================================================
% Step 1.2: NERNST-PLANCK - bound for \dert\cl
%====================================================================
\underline{Step~1.2 -- estimates for $\dert\cl$:}~~
We abbreviate $B:=\Lp[I]{2}{;\spaceTC}$. By involving equations~\eqref{eq:transportWeak-fpi}, we obtain the identity
\begin{align*}
 &     \theta\norm{\dert\cl}{\Lp[I]{2}{;\spaceTC^\ast}}:=\sup_{\norm{\test}{B}\leq1} \dualp[]{\theta\dert\cl}{\test}_{B^\ast\times B} \\
 &=    \sup_{\norm{\test}{B}\leq1} \sqbrac{ -\scp[]{\Dl\grad\cl}{\grad\test}_{\OmegaT} + \scp[]{{\cl}[\fieldF+\coeffEL\fieldEL]}{\grad\test}_{\OmegaT} + \scp[]{\theta\Rl(\Cl)}{\test}_{\OmegaT} + \scp[]{\gl}{\test}_{\GammaT} } \\
 &=: I.1 ~+~ I.2 ~+~ I.3 ~+~ I.4~.
\end{align*}
For $I.1$ and $I.3$, we arrive with \Holder, \ref{Assump:Ellip}, \ref{Assump:Reaction}, and \cref{lemma:bounded} at
\begin{align*}
        I.1 + I.3
  &\leq C_D\norm{\grad\cl}{\Lp[\OmegaT]{2}{}} \sqbrac{\sup_{\norm{\test}{B}\leq1} \norm{\grad\test}{\Lp[\OmegaT]{2}{}}} +C_{\Rl}\sum_l \norm{\cl}{\Lp[\OmegaT]{2}{}} \sqbrac{\sup_{\norm{\test}{B}\leq1}\norm{\test}{\Lp[\OmegaT]{2}{}}} \\
  &\leq (C_D+C_{\Rl}) \sum_l \norm{\cl}{B}
    ~\leq~ (C_D+\max_l C_{\Rl})C_0(\timeEnd) =: C_1(\timeEnd) .
\end{align*}
The integral~$I.2$, we bound with \Holder\ and \cref{lemma:bounded} by
\begin{align*}
        I.2
 &\leq \norm{{\cl}[\fieldF+\coeffEL\fieldEL]}{\Lp[\OmegaT]{2}{}} \sqbrac{ \sup_{\norm{\test}{B}\leq1}\norm{\grad\test}{\Lp[\OmegaT]{2}{}}} \\
 &\leq \norm{\cl}{\Lp[\OmegaT]{\infty}{}} \norm{\fieldF+\coeffEL\fieldEL}{\Lp[\OmegaT]{2}{}}
 ~\leq~ C_M(\timeEnd) \norm{\fieldF+\coeffEL\fieldEL}{\Lp[\OmegaT]{2}{}} ~.
\end{align*}
For $I.4$, we immediately get with \Holder\ and \cref{lemma:interpol-boundary}
\begin{align*}
      I.4
 \leq \norm{\gl}{\Lp[\GammaT]{2}{}} \sqbrac{ \sup_{\norm{\test}{B}\leq1} \norm{\test}{\Lp[\GammaT]{2}{}} }
 \leq \norm{\gl}{\Lp[\GammaT]{2}{}} \sqbrac{ C\sup_{\norm{\test}{B}\leq1} \norm{\test}{\Lp[\OmegaT]{2}{}} }
 \leq C\norm{\gl}{\Lp[\GammaT]{2}{}}~.
\end{align*}
Thus, by combining the estimates for $I.1$ -- $I.4$, we have shown
\begin{align}\label{eq:apriori-np-time}
  \norm{\dert\cl}{\Lp[I]{2}{;\spaceTC^\ast}} \leq
  C_1(\timeEnd) + C_M(\timeEnd) \norm{\fieldF+\coeffEL\fieldEL}{\Lp[\OmegaT]{2}{}} + C\norm{\gl}{\Lp[\GammaT]{2}{}}.
\end{align}
%
%====================================================================
% Step 1.3: NERNST-PLANCK - final a priori bound
%====================================================================
\underline{Step~1.3 -- a priori estimates for $\cl$:}~~
We now put \cref{lemma:energy}, \cref{lemma:bounded}, and the estimates~\eqref{eq:apriori-np-time} together.
In anticipation of estimates~\eqref{eq:apriori-gauss} and \eqref{eq:apriori-darcy}, we obtain the desired a~priori bound
\begin{align*}
 &     \sum_l\sqbrac{ \norm{\cl}{\Lp[I]{\infty}{;\Lp{2}{}}}+\norm{\cl}{\Lp[I]{2}{;\spaceTC}}+\norm{\cl}{\Hk[I]{1}{;\spaceTC^\ast}}+\norm{\cl}{\Lp[\OmegaT]{\infty}{}} } \\
 &\leq C_0(\timeEnd) + C_M(\timeEnd) + C_1(\timeEnd) + \max(1,\coeffEL)(C_e+C_f)C_M(\timeEnd) + C\norm{\gl}{\Lp[\GammaT]{2}{}} .
\end{align*}
%
%====================================================================
% Step 2.1: GAUSS - bound for \grad\cdot\fieldEL
%====================================================================
\underline{Step~2.1 -- estimate for $\grad\cdot\fieldEL$:}~~
We test equation~\eqref{eq:gaussWeakDIV-fpi} with $\test=\grad\cdot\fieldEL$. Thereby, we directly get with \Young\
\begin{align*}
  \norm{\grad\cdot\fieldEL}{\Lp{2}{}}^2 \leq \norm{\rho_b}{\Lp{2}{}}^2 + \theta\max_l\abs{\zl} \sum_l\norm{\clold}{\Lp{2}{}}^2~.
\end{align*}
Since this estimate holds uniformly in time, we take the supremum over $t\in[0,\timeEnd]$ and come to
\begin{align*}
  \norm{\grad\cdot\fieldEL}{\Lp[I]{\infty}{;\Lp{2}{}}}^2 \leq \norm{\rho_b}{\Lp[I]{\infty}{;\Lp{2}{}}}^2 + \theta\max_l\abs{\zl} \sum_l\norm{\clold}{\Lp[I]{\infty}{;\Lp{2}{}}}^2~.
\end{align*}
Hence, together with \cref{lemma:energy} we finally have (assume $\norm{\clold}{\Lp[I]{\infty}{;\Lp{2}{}}}\leq R$ in case $\clold\neq\cl$)
\begin{align*}
       \norm{\grad\cdot\fieldEL}{\Lp[I]{\infty}{;\Lp{2}{}}}
 &\leq \begin{cases} \norm{\rho_b}{\Lp[I]{\infty}{;\Lp{2}{}}} + \theta\max_l\abs{\zl} R      & \text{ if } \cl\neq\clold,\\[2.0mm]
                     \norm{\rho_b}{\Lp[I]{\infty}{;\Lp{2}{}}} + \theta\max_l\abs{\zl} C_0    & \text{ if } \cl=\clold.
       \end{cases}\\[2.0mm]
 &=:C_{1,e}(\timeEnd,R)~.
\end{align*}
%
%====================================================================
% Step 2.2: GAUSS - bound for \potEL
%====================================================================
\underline{Step~2.2 -- estimate for $\potEL$:}~~
Next, we test equation~\eqref{eq:gaussWeak-fpi} with $\Test\in\spaceTF$. Due to \cite[Chapter~7.2]{Quarteroni-book}, we can choose $\Test$ such that
$\grad\cdot\Test=\potEL$ and $\norm{\Test}{\Hkdiv{1}{}}\leq K\norm{\potEL}{\spaceTP}$ holds. This yields with \ref{Assump:Ellip} and \Young\ 	
\begin{align*}
   \norm{\potEL}{\Lp{2}{}}^2 \leq \frac{K^2}{\permitELscalar\alpha_D}\norm{\fieldEL}{\Lp{2}{}}^2
   ~~~~\Hence~
   \norm{\potEL}{\Lp[I]{\infty}{;\Lp{2}{}}}^2 \leq \frac{K^2}{\permitELscalar\alpha_D}\norm{\fieldEL}{\Lp[I]{\infty}{;\Lp{2}{}}}^2~.
\end{align*}
%
%====================================================================
% Step 2.3: GAUSS - bound for \fieldEL
%====================================================================
\underline{Step~2.3 -- estimate for $\fieldEL$:}~~
Due to $\fieldEL\in\spaceSEF$ and \ref{Assump:BoundData}, we test equation~\eqref{eq:gaussWeak-fpi} with $\Test=\fieldEL-\vecsigma\in\spaceTF$. In addition, we test equation~\eqref{eq:gaussWeakDIV-fpi}  with $\test=\potEL$.
% This leads us to
% \begin{align*}
%  \scp{\permitEL^{-1}\fieldEL}{\fieldEL}_\Omega &= \scp{\potEL}{\grad\cdot\fieldEL}_\Omega + \scp{\permitEL^{-1}\fieldEL}{\vecsigma}_\Omega - \scp{\potEL}{\grad\cdot\vecsigma}_\Omega,   \\
%  \scp{\grad\cdot\fieldEL}{\potEL}_{\Omega}     &= \scp{\rho_b + \chargeELold}{\potEL}_\Omega .
% \end{align*}
By adding these equations, we get with \ref{Assump:Ellip} and \Young\
\begin{align*}
 &     \frac{1}{\permitELscalar C_D}\norm{\fieldEL}{\Lp{2}{}}^2
    ~= \scp{\permitEL^{-1}\fieldEL}{\fieldEL}_\Omega
 ~\leq \scp{\permitEL^{-1}\fieldEL}{\vecsigma}_\Omega - \scp{\potEL}{\grad\cdot\vecsigma}_\Omega+ \scp{\rho_b + \chargeELold}{\potEL}_\Omega \\
 &\leq \frac{\delta_1}{\permitELscalar\alpha_D}\norm{\fieldEL}{\Lp{2}{}}^2 + \frac{1}{4\delta_1}\norm{\vecsigma}{\Lp{2}{}}^2 + \delta_2\norm{\potEL}{\Lp{2}{}}^2 + \frac{1}{4\delta_2} \sqbrac{\norm{\grad\cdot\vecsigma}{\Lp{2}{}}^2 + \norm{\rho_b + \chargeELold}{\Lp{2}{}}^2}.
\end{align*}
Thus, we arrive with a suitable choice of $\delta_1$, $\delta_2$, the above estimate for $\potEL$, by taking the supremum over time,
and with \cref{lemma:bounded} at (we assume $\cl=\clold$ and we skip the uncritical case $\cl\neq\clold$)
\begin{align*}
       \norm{\fieldEL}{\Lp[I]{\infty}{;\Lp{2}{}}}
 \leq  \kappa\sqbrac{\norm{\vecsigma}{\Lp[I]{\infty}{;\Hkdiv{1}{}}} + \norm{\rho_b}{\Lp[I]{\infty}{;\Lp{2}{}}} +\theta\max_l\abs{\zl}C_0 }
 ~:= C_{2,e}(\timeEnd).
\end{align*}
%
%====================================================================
% Step 2.4: GAUSS - final a priori bound
%====================================================================
\underline{Step~2.4 -- a priori estimate for $(\fieldEL,\potEL)$:}~~
Collecting the preceding inequalities for $\grad\cdot\fieldEL$, $\fieldEL$, and $\potEL$ shows
\begin{align}\label{eq:apriori-gauss}
      \norm{\fieldEL}{\Lp[I]{\infty}{;\Hkdiv{1}{}}} + \norm{\potEL}{\Lp[I]{\infty}{;\Lp{2}{}}}
 \leq C_{1,e} + C_{2,e} + \frac{K}{\sqrt{\permitELscalar\alpha_D}} C_{2,e}
 ~=:  C_e(\timeEnd).
\end{align}
%
%====================================================================
% Step 3.1: DARCY - bound for \grad\cdot\fieldF
%====================================================================
\underline{Step~3.1 -- estimate for $\grad\cdot\fieldF$:}~~
We test equation~\eqref{eq:darcyWeakDIV-fpi} with $\test=\grad\cdot\fieldF$ and immediately obtain $\norm{\grad\cdot\fieldF}{\Lp{2}{}}^2 =0$ and thus $\norm{\grad\cdot\fieldF}{\Lp[I]{\infty}{;\Lp{2}{}}} =0$.\\[2.0mm]
%
%====================================================================
% Step 3.2: DARCY - bound for \press
%====================================================================
\underline{Step~3.2 -- estimate for $\press$:}~~
Next, we test equation~\eqref{eq:darcyWeak-fpi} with $\Test\in\spaceTF$. According to \cite[Chapter~7.2]{Quarteroni-book}, we find a $\Test$ such that
$\grad\cdot\Test=\press$ and $\norm{\Test}{\Hkdiv{1}{}}\leq K\norm{\press}{\spaceTP}$ holds. This leads us  with \ref{Assump:Ellip}, \Young, \cref{lemma:bounded}, and \eqref{eq:apriori-gauss} to
(we assume $\cl=\clold$ and we skip the uncritical case $\cl\neq\clold$)
\begin{align*}
   \norm{\press}{\Lp{2}{}}^2  &\leq \delta K \norm{\press}{\Lp{2}{}}^2 + \frac{\mu C_K}{2\delta}\norm{\fieldF}{\Lp{2}{}}^2 +\frac{1}{2\delta\permitELscalar\alpha_D}\norm{\chargeELold\fieldEL}{\Lp{2}{}}^2 \\
                              &\leq \delta K \norm{\press}{\Lp{2}{}}^2 + \frac{\mu C_K}{2\delta}\norm{\fieldF}{\Lp{2}{}}^2 +\frac{\theta\max_l\abs{\zl}}{2\delta\permitELscalar\alpha_D} C_e^2 C_M^2~.
\end{align*}
A suitable choice of $\delta>0$ immediately shows
\begin{align*}
 \norm{\press}{\Lp[I]{\infty}{;\Lp{2}{}}} &\leq 2\mu K C_K \norm{\fieldF}{\Lp[I]{\infty}{;\Lp{2}{}}} + \underbrace{\frac{2K\theta\max_l\abs{\zl}}{\permitELscalar\alpha_D} C_e(\timeEnd) C_M(\timeEnd)}_{:=C_{1,f}}  ~.
\end{align*}
%
%====================================================================
% Step 3.3: DARCY - bound for \fieldF
%====================================================================
\underline{Step~3.3 -- estimate for $\fieldF$:}~~
We test equation~\eqref{eq:darcyWeakDIV-fpi} with $\test=\press$ and equation~\eqref{eq:darcyWeak-error} with the test function~$\Test=\fieldF-\vecf$. Here, we take $\vecf$ according to \ref{Assump:BoundData},
which ensures $\Test\in\spaceTF$.
% Thereby, we come to
% \begin{align*}
%   & \scp{\permeabH^{-1}\fieldF}{\fieldF}_\Omega = \scp{\mu^{-1}\press}{\grad\cdot\sqbrac{\fieldF-\vecf}}_\Omega + \scp{\permeabH^{-1}\fieldF}{\vecf}_\Omega + \scp{\coeffForceEL\forceELold}{\fieldF-\vecf}_\Omega,   \\
%   & \scp{\grad\cdot\fieldF}{\press}_\Omega = 0 ~.
% \end{align*}
Furthermore, adding these equations, yields with \ref{Assump:Ellip} and \Young, \cref{lemma:bounded}, and \eqref{eq:apriori-gauss} (again, we assume $\cl=\clold$ and we skip the uncritical case $\cl\neq\clold$)
\begin{align*}
        \alpha_K \norm{\fieldF}{\Lp{2}{}}^2
  &\leq %\scp[]{\permeabH^{-1}\fieldF}{\fieldF}_\Omega =
        -\scp{\mu^{-1}\press}{\grad\cdot\vecf}_\Omega + \scp{\permeabH^{-1}\fieldF}{\vecf}_\Omega + \scp{\coeffForceEL\forceELold}{\fieldF-\vecf}_\Omega \\
  &\leq \frac{\delta_1}{2}\norm{\press}{\Lp{2}{}}^2 + \frac{1}{2\mu^2\delta_1}\norm{\grad\cdot\vecf}{\Lp{2}{}}^2 + \delta_2\norm{\fieldF}{\Lp{2}{}}^2 + \brac{\frac{C_K^2}{2\delta_2}+\frac{1}{2}} \norm{\vecf}{\Lp{2}{}}^2 \\
  &~~~~ +\sqbrac{\frac{1}{2\permitELscalar\alpha_D\mu\delta_2} + \frac{1}{2}} \norm{\forceELold}{\Lp{2}{}}^2 \\
  &\leq \frac{\delta_1}{2}\norm{\press}{\Lp{2}{}}^2 + \frac{1}{2\mu^2\delta_1}\norm{\grad\cdot\vecf}{\Lp{2}{}}^2 + \delta_2\norm{\fieldF}{\Lp{2}{}}^2 + \brac{\frac{C_K^2}{2\delta_2}+\frac{1}{2}} \norm{\vecf}{\Lp{2}{}}^2 \\
  &~~~~ +\theta\max_l\abs{\zl}\sqbrac{ \frac{1}{2\permitELscalar\alpha_D\mu\delta_2} + \frac{1}{2}}C_e^2 C_M^2 ~.
\end{align*}
We now insert the estimate for $\press$ and we choose $\delta_1$ and $\delta_2$ appropriately. Thereby, we directly arrive with taking the supremum over time at
% \begin{align*}
%         \norm{\fieldF}{\Lp{2}{}}^2
%   &\leq \frac{4K C_K}{\mu\alpha^2_K}\norm{\grad\cdot\vecf}{\Lp{2}{}}^2 + \brac{\frac{4C_K^2}{\alpha^2_K}+\frac{2}{\alpha_K}} \norm{\vecf}{\Lp{2}{}}^2 \\
%   &~~~~ +\theta\max_l\abs{\zl}\sqbrac{\frac{8K}{\permitELscalar\alpha_D\alpha_K}+\frac{1}{\permitELscalar\alpha_D\alpha^2_K\mu} + \frac{2}{\alpha_K}}C_e^2 C_M^2 \\
%   &\leq \brac{\frac{4K C_K}{\mu\alpha^2_K}+\frac{4C_K^2}{\alpha^2_K}+\frac{2}{\alpha_K}} \norm{\vecf}{\Hkdiv{1}{}}^2 \\
%   &~~~~ +\theta\max_l\abs{\zl}\sqbrac{\frac{8K}{\permitELscalar\alpha_D\alpha_K}+\frac{1}{\permitELscalar\alpha_D\alpha^2_K\mu} + \frac{2}{\alpha_K}}C_e^2 C_M^2.
% \end{align*}
% Thus, we know that
\begin{align*}
        \norm{\fieldF}{\Lp[I]{\infty}{;\Lp{2}{}}}
  &\leq \brac{\frac{4K C_K}{\mu\alpha^2_K}+\frac{2C_K}{\alpha_k}+\frac{2}{\alpha_K}} \norm{\vecf}{\Hkdiv{1}{}}^2 \\
  &~~~~ +\theta\max_l\abs{\zl}\sqbrac{\frac{8K}{\permitELscalar\alpha_D\alpha_K}+\frac{1}{\permitELscalar\alpha_D\alpha_K\mu} + \frac{2}{\alpha_K}}C_e C_M
  ~=: C_{2,f}(\timeEnd).
\end{align*}
%
%====================================================================
% Step 3.4: DARCY - final a priori bound
%====================================================================
\underline{Step~3.4 -- a priori estimate for $(\fieldF,\press)$:}~~
Combing the preceding estimates for $\grad\cdot\fieldF$, $\fieldF$, $\press$ shows
\begin{align}\label{eq:apriori-darcy}
     \norm{\fieldF}{\Lp[I]{\infty}{;\Hkdiv{1}{}}} + \norm{\press}{\Lp[I]{\infty}{;\Lp{2}{}}}
 \leq C_{2,f} + 2\mu K C_KC_{2,f} +  C_{1,f}  ~=: C_f(\timeEnd)~.
\end{align}

\end{proof}
\begin{remark}
\textnormal{
The proof of \cref{thm:aprioriBounds} is valid in arbitrary space dimensions, i.e., for $\Omega\subset\setR^n$ with $n\geq2$.
However, in \ref{Assump:Geom} we restrict ourselves to $n\leq3$, as we use in the proof of \cref{thm:existence} compact embeddings
of Aubin-Lions-type, which are valid only for $n\leq3$.
} \hfill$\square$
\end{remark}
%
%
%%%%%%%%%%%%%%%%%%%%%%%%%%%%%%%%%%%%%%%%%%
%  SUBSECTION existence of a fixed point %
%%%%%%%%%%%%%%%%%%%%%%%%%%%%%%%%%%%%%%%%%%
\subsection{Existence of a fixed point}\label{subsec:fpi-point}
In this section, we prove the  existence of global weak solutions of the \dpnp. Our proof is based on the following fixed point theorem, see \cite[Corollary 9.6]{Zeidler1-book}.
%
%===============================
%===============================
% THEOREM: FIXED POINT THEOREM
%===============================
%===============================
\begin{thm}\label{thm:fpi-theorem}
Let $\fixOp: \fixSet\subset X\rightarrow \fixSet$ be continuous, where $\fixSet$ is a nonempty, compact, and convex set in a locally convex space $X$. Then, $\fixOp$ has a fixed point.
\end{thm}
A Banach space~$X$ equipped with the $\weakstar$-topology is a locally convex space $(X,\weakstar)$. Hence, the above fixed point theorem is tailored for Banach spaces,
which carry the $\weakstar$-topology.  In our case, the $\weakstar$-topology is the natural choice for the following three reasons:
\par
Firstly, the a~priori estimates from \cref{subsec:aprioriBounds} are equivalent to $\weakstar$-compactness.
Secondly, the solution space for $\cl$ includes $\Lp[\OmegaT]{\infty}{}$, which is not reflexive. Hence, the $\weakstar$-topology differs
from the $\weak$-topology. Thirdly, when using the $\weakstar$-topology, we can reuse the a~priori estimates from \cref{subsec:aprioriBounds}
for the $\weakstar$-continuity of the fixed point operator.
\par
In summary, we can exaggeratedly state that in the $\weakstar$-topology, the compactness of $\fixSet$ and the continuity of $\fixOp$ is already contained in the a~priori estimates.
However, this is valid, only if the predual of the solution space is separable. In this case, the set-based topological terms and the sequences-based ones coincide.
This enables us to prove the continuity of the operator with $\weakstar$-convergent sequences, instead of investigating preimages of $\weakstar$-open sets.
%
%============================
%============================
% THEOREM: EXISTENCE
%============================
%============================
\begin{thm}\label{thm:existence}
Let \ref{Assump:Geom}--\ref{Assump:back-charge} be valid. Then, there exists a solution $\brac{\fieldEL,\potEL,\fieldF,\press,\Cl}\in \setR^{4+2n}$ of equations~\eqref{eq:Model-1a}--\eqref{eq:Model-1j}
according to \cref{def:weaksolution}.
\end{thm}
\begin{proof}
For ease of readability, we split the proof into several steps\\[2.0mm]
%
%====================================================================
% Step 1 - construct the space X
%====================================================================
\underline{Step~1 -- the space $X$:}~~
First of all, we repeat the definition of the space~$X$ from \cref{def:FixOp}
\begin{align*}
 X :=\sqbrac{\fspace{L^\infty}{I}{;\Lp{2}{}} \cap \spaceSC \cap \Lp[\OmegaT]{\infty}{}}^2 .
\end{align*}
Furthermore, we equip $X$ with the norm
\begin{align*}
 \norm{\cdot}{X} := \norm{\cdot}{\fspace{L^\infty}{I}{;\Lp{2}{}}} + \norm{\cdot}{\fspace{L^2}{I}{;\Hk{1}{}}} + \norm{\cdot}{\fspace{H^¹}{I}{;\Hk{1}{}}^\ast} + \norm{\cdot}{\Lp[\OmegaT]{\infty}{}} .
\end{align*}
Thus, $(X,\norm{\cdot}{X})$ is a Banach space. However, we henceforth consider the locally convex space $(X,\weakstar)$  and all topological terms refer to the $\weakstar$-topology.
Furthermore, the predual $X_0$ of $X$ can be written according to \cite[Chapter I, IV]{Gajewski-book}) as
\begin{align*}
X_0 := \sqbrac{\fspace{L^1}{I}{;\Lp{2}{}} + \fspace{L^2}{I}{;\Hk{1}{}^\ast} +\fspace{H^1}{I}{;\Hk{1}{}}^\ast + \Lp[\OmegaT]{1}{}}^2 .
\end{align*}
Hence, $X_0$ is a separable Banach space with dual $X$ and the topological terms for $(X,\weakstar)$ based on sets are equivalent with those based on sequences,~cf.~\cite{Zeidler1-book,Gajewski-book}.
In particular, the notion of $\weakstar$-continuous/compact is equal to sequentially $\weakstar$-continuous/compact.\\[2.0mm]
%
%====================================================================
% Step 2 - the set K
%====================================================================
\underline{Step~2 -- the set $\fixSet$:}~~
For $R>0$, we introduce the set $\fixSet$ as a ball of radius $R$ in $X$, i.e.,
\begin{align*}
\fixSet := \cbrac{ \vecv=(v_1,v_2)\in X: ~~\norm{\vecv}{X}\leq R } \subset X~.
\end{align*}
$\fixSet$ is nonempty, convex, and $\weakstar$-compact due to Banach-Alaoglu-Bourbaki~theorem, cf. \cite[Theorem 1.7]{Roubicek-book}.\\[2.0mm]
%
%====================================================================
% Step 3 - the operator F
%====================================================================
\underline{Step~3 -- the operator $\fixOp$:}~~
We consider the operator~$\fixOp$, which was already introduced in \cref{def:FixOp}. This operator is a well-defined operator due to \cref{lemma:wellDef-fixOp}.\\[2.0mm]
%
%====================================================================
% Step 4 - the self mapping property
%====================================================================
\underline{Step~4 -- self mapping property $\fixOp(\fixSet)\subset\fixSet$:}~~
Let $\Clold\in\fixSet$. The definition of the set~$\fixSet$ and the definition of the norm~$\norm{\cdot}{X}$ ensure that we have for $l=1,2$
\begin{align*}
\norm{\clold}{\fspace{L^\infty}{I}{;\Lp{2}{}}} + \norm{\clold}{\fspace{L^2}{I}{;\Hk{1}{}}} + \norm{\clold}{\fspace{H^1}{I}{;\Hk{1}{}^\ast}} + \norm{\clold}{\Lp[\OmegaT]{\infty}{}} \leq R~.
\end{align*}
With this information, we return to the a~priori estimates in \cref{thm:aprioriBounds}. Carefully reading through the proof of \cref{thm:aprioriBounds} reveals in detail how the constants
of the a~priori bounds are defined. More precisely, this shows that
\begin{align*}
 \norm{\Cl}{X} \text{ is bounded in terms of } \begin{cases} \text{the data, the radius~$R$ and the end time~$\timeEnd$} & \text{if } \Cl\neq\Clold, \\
                                                                \text{the data and the end time~$\timeEnd$}                 & \text{if } \Cl=\Clold.
                                                   \end{cases}
\end{align*}
In both cases, the constants of the a~priori estimate are partially independent of the end time~$\timeEnd$. This means, we can split the constants into
\begin{align*}
 C = C(\timeEnd) + C(\timeEnd,R) + C_d~.
\end{align*}
We now choose the radius~$R:=2C(\timeEnd) + 2C_d$~. In the remaining part~$C(\timeEnd,R)$, we assume a sufficiently small end time~$\timeEnd<<1$ such that we have $C(\timeEnd,R) \leq C(\timeEnd)+C_d$~.
This proves
\begin{align*}
 \norm{\Cl}{X} ~\leq~ C ~=~  C(\timeEnd) + C(\timeEnd,R) + C_d ~\leq ~ 2C(\timeEnd) + 2C_d = R ~.
\end{align*}
Thus, we have $\fixOp(\fixSet)\subset\fixSet$. However, we note that due to the assumption $\timeEnd<<1$, we are restricted to generally small time intervals $I=[0,\timeEnd]$.\\[2.0mm]
%
%====================================================================
% Step 5 - the weakstar continuity
%====================================================================
\underline{Step~5 -- $\weakstar$-continuity of $\fixOp$}\\
Subsequently, we use the already mentioned equivalence between $\weakstar$-continuous and sequentially $\weakstar$-continuous. This means,
we show the $\weakstar$-continuity with the criterion based on sequences. For that purpose, we consider a sequence $(\Clold)_k\subset\fixSet$, for which we assume that
$\Clold_k\overset{*~}\rightharpoonup\Clold$ in $X$. As $\fixOp(\Clold_k)=\Cl_k$ is the solution of \eqref{eq:transportWeak-fpi}, we know together with $\Clold_k\in\fixSet$
and the just established self mapping property that
\begin{align}\label{eq:existence-conv-1}
 \Clold_k ~\overset{*~}\rightharpoonup~ \Clold ~~\text{ in } X \text{ with }~ \Clold\in\fixSet \qquad\text{and}\qquad \norm{\Cl_k}{X}=\norm{\fixOp(\Clold_k)}{X}\leq R~.
\end{align}
Consequently, $(\Cl_k)_k$ is a uniformly bounded sequence and a subsequence, denoted again by $(\Cl_k)_k$,  $\weakstar$-converges to a unique limit $\Cl\in\fixSet$, i.e., it holds
\begin{align}\label{eq:existence-conv-2}
 \Cl_k ~\overset{*~}\rightharpoonup~ \Cl ~~\text{ in } X \text{ with }~ \Cl\in\fixSet~.
\end{align}
The fixed point operator~$\fixOp$ is $\weakstar$-continuous, if and only if $\Cl$ solves the \enquote{limit} PDE \eqref{eq:transportWeak-fpi}, which is generated by $\Clold$. In this case, we have
\begin{align}\label{eq:existence-defWeakstarContinuous}
 \Clold_k~\overset{*~}\rightharpoonup~ \Clold
 ~~\Hence~~
 \Cl_k= \fixOp(\Clold_k) ~\overset{*~}\rightharpoonup~ \fixOp(\Clold)=\Cl \quad \text{ in } X~.
\end{align}
Thus, it remains to show that $\Cl$ solves the \enquote{limit} PDE \eqref{eq:transportWeak-fpi}, which is generated by $\Clold$. To show this, we return to \cref{def:FixOp}
and we subtract the equations, which are generated by $\Clold_k$ from the equations, which are generated by $\Clold$ and we integrate in time. Thereby, we obtain the error equations\\[2.0mm]
%===================
%% error equations
%==================
\begin{subequations}
\underline{Gauss's law:}\vspace{-2mm}
\begin{align}
 \scp{\permitEL^{-1}(\fieldEL_k-\fieldEL)}{\Test}_{\OmegaT} &= \scp{\potEL_k-\potEL}{\grad\cdot\Test}_{\OmegaT}, \label{eq:gaussWeak-weakstarError} \\
   \scp{\grad\cdot(\fieldEL_k-\fieldEL)}{\test}_{\OmegaT} &= \theta \scp{\chargeELold[,k]-\chargeELold}{\test}_{\OmegaT}~. \label{eq:gaussWeakDIV-weakstarError}
\end{align}
\underline{Darcy's law:}\vspace{-2mm}
\begin{align}
 \scp{\permeabH^{-1}(\fieldF_k-\fieldF)}{\Test}_{\OmegaT} &= \scp{\mu^{-1}(\press_k-\press)}{\grad\cdot\Test}_{\OmegaT} + \theta\mu^{-1}\scp{\chargeELold[,k]\permitEL^{-1}\fieldEL_k}{\Test}_{\OmegaT} \nonumber\\
                                                          &~~~~-\theta\mu^{-1}\scp{\chargeELold\permitEL^{-1}\fieldEL}{\Test}_{\OmegaT}, \label{eq:darcyWeak-weakstarError}\\
     \scp{\grad\cdot(\fieldF_k-\fieldF)}{\test}_{\OmegaT} &= 0~. \label{eq:darcyWeakDIV-weakstarError} \\[6.0mm] \nonumber
\end{align}
\underline{Nernst--Planck equations:}\vspace{-2mm}
\begin{align}
   & \dualp{\dert(\cl[l,k] - \cl)}{\test}_{\Lp[I]{2}{;\Hk{1}{}^\ast}\times\Lp[I]{2}{;\Hk{1}{}}} + \scp{\Dl\grad(\cl[l,k]-\cl)}{\grad\test}_{\OmegaT} \nonumber\\[2.0mm]
 &  ~~-\scp{ {\cl[l,k]}\sqbrac{\fieldF_k +\coeffEL\fieldEL_k}-{\cl}\sqbrac{\fieldF+\coeffEL\fieldEL} }{\grad\test}_{\OmegaT} \nonumber\\[2.0mm]
 &=  \theta\scp{\Rl(\Cl_k)-\Rl(\Cl)}{\test}_{\OmegaT} ~. \label{eq:transportWeak-weakstarError}
\end{align}
\end{subequations}
%===================
%% Aubin-Lions
%==================
We note that \eqref{eq:existence-conv-1} and Aubin-Lions Lemma, cf.~\cite[Lemma 7.7]{Roubicek-book}, imply the norm-convergences
\begin{align}\label{eq:existence-AubinLions}
 \clold[l,k] \rightarrow \clold ~\text{ in } \Lp[\OmegaT]{2}{}  \qquad\text{and}\qquad \clold[l,k] \rightarrow \clold ~\text{ in } \Lp[I]{2}{;\Lp{3}{}}~.
\end{align}
%===================
%% Error Gauss
%==================
Hence, we obtain for $(\fieldEL,\potEL)$  analogously to the proof of \cref{thm:unique} or \cref{thm:aprioriBounds}, the norm-convergence
$\norm{\potEL_k-\potEL}{\Lp[\OmegaT]{2}{}} + \norm{\fieldEL_k-\fieldEL}{\Lp[I]{2}{;\Hkdiv{1}{}}} \leq C \sum_l\norm{\clold[l,k]-\clold}{\Lp[\OmegaT]{2}{}} ~~\rightarrow 0~.$
\par
%===================
%% Error Darcy
%==================
Furthermore, for $(\fieldF,\press)$, we get analogously to the proof of \cref{thm:unique} or \cref{thm:aprioriBounds}
\begin{align*}
  \norm{\press_k-\press}{\Lp[\OmegaT]{2}{}} + \norm{\fieldF_k-\fieldF}{\Lp[I]{2}{;\Hkdiv{1}{}}} \leq C \norm{\chargeELold[,k]\permitEL^{-1}\fieldEL_k-\chargeELold\permitEL^{-1}\fieldEL}{\Lp[\OmegaT]{2}{}}~.
\end{align*}
Applying \ref{Assump:Ellip}, \Holder, \cref{lemma:regularity-gausslaw}, and \eqref{eq:existence-AubinLions} yields for the right hand side
\begin{align*}
  &     C\norm{\chargeELold[,k]\permitEL^{-1}\fieldEL_k-\chargeELold\permitEL^{-1}\fieldEL}{\Lp[\OmegaT]{2}{}} \\
  &\leq C\norm{\chargeELold[,k]\permitEL^{-1}(\fieldEL_k-\fieldEL)}{\Lp[\OmegaT]{2}{}} +  C\norm{(\chargeELold[,k]-\chargeELold)\permitEL^{-1}\fieldEL}{\Lp[\OmegaT]{2}{}} \\
 % &\leq C\norm{\chargeELold[,k]}{\Lp[\OmegaT]{\infty}{}} \norm{\fieldEL_k-\fieldEL}{\Lp[\OmegaT]{2}{}} + C\norm{\fieldEL}{\Lp[I]{2}{;\Lp{6}{}}} \norm{\chargeELold[,k]-\chargeELold}{\Lp[I]{2}{;\Lp{3}{}}} \\
  &\leq CR\norm{\fieldEL_k-\fieldEL}{\Lp[\OmegaT]{2}{}} + C\norm{\fieldEL}{\Lp[I]{2}{;\Hk{1}{}}} \norm{\chargeELold[,k]-\chargeELold}{\Lp[I]{2}{;\Lp{3}{}}}
   ~~\rightarrow 0~.
\end{align*}
This shows the norm-convergence $\norm{\press_k-\press}{\Lp[\OmegaT]{2}{}} + \norm{\fieldF_k-\fieldF}{\Lp[I]{2}{;\Hkdiv{1}{}}} ~~\rightarrow 0~.$
\par
%===================
%% error NP
%==================
As to the convergence for $\Cl_k$, we begin with the time integrals and the diffusion integrals. From \eqref{eq:existence-conv-2} follows that $\dert\cl[l,k]$ resp. $\grad\cl[l,k]$
$\weakstar$-converge%
\footnote{In both cases $\weakstar$-convergence is equivalent to $\weak$-convergence as the involved spaces are reflexive.}
towards $\dert\cl$ in $\Lp[I]{2}{;\Hk{1}{}^\ast}$ resp. $\grad\cl$ in $\Lp[\OmegaT]{2}{}$. Thus, we have
\begin{align*}
  \dualp{\dert(\cl[l,k] - \cl)}{\test}_{\Lp[I]{2}{;\Hk{1}{}^\ast}\times\Lp[I]{2}{;\Hk{1}{}}} + \scp{\Dl\grad(\cl[l,k]-\cl)}{\grad\test}_{\OmegaT}  \rightarrow 0.
\end{align*}
For the convection integrals and the electric~drift integrals, we obtain
\begin{align*}%\label{eq:existence-flow&drift}
 &   \scp{ {\cl[l,k]}\sqbrac{\fieldF_k +\coeffEL\fieldEL_k}-{\cl}\sqbrac{\fieldF+\coeffEL\fieldEL} }{\grad\test}_{\OmegaT} \nonumber\\
 &=  \scp{ ({\cl[l,k]}-{\cl})\sqbrac{\fieldF +\coeffEL\fieldEL} }{\grad\test}_{\OmegaT} \nonumber\\
 &   ~~~~+\scp{ {\cl[l,k]}\sqbrac{(\fieldF_k-\fieldF) +\coeffEL(\fieldEL_k-\fieldEL)} }{\grad\test}_{\OmegaT}
 ~~=: I.1 + I.2~.
\end{align*}
Concerning $I.1$, we note that $\sqbrac{\fieldF +\coeffEL\fieldEL}\grad\test \in \Lp[\OmegaT]{1}{}$ and $\cl[l,k]-\cl\in\Lp[\OmegaT]{\infty}{}$ with $\Lp[\OmegaT]{\infty}{}=\Lp[\OmegaT]{1}{}^\ast$.
Thereby, we arrive with \eqref{eq:existence-conv-2} at
\begin{flalign*}
 I.1 %&=\Intdt{\Intdx{({\cl[l,k]}-{\cl})\sqbrac{\fieldF +\coeffEL\fieldEL}\cdot\grad\test ~}}  \\
     &= \dualp{({\cl[l,k]}-{\cl})}{\sqbrac{\fieldF +\coeffEL\fieldEL}\grad\test}_{\Lp[\OmegaT]{1}{}^\ast\times\Lp[\OmegaT]{1}{}} \rightarrow 0 .
\end{flalign*}
Concerning $I.2$, we come with \Holder\ and the $\Lp[\OmegaT]{2}{}$-convergence of $\fieldEL_k$, $\fieldF_k$ to
\begin{align*}
 I.2 %&\leq \norm{\cl[l,k]}{\Lp[\OmegaT]{\infty}{}} \brac{\norm{\fieldF_k-\fieldF}{\Lp[\OmegaT]{2}{}}+\coeffEL\norm{\fieldEL_k-\fieldEL}{\Lp[\OmegaT]{2}{}}} \norm{\grad\test}{\Lp[\OmegaT]{2}{}} \\
     &\leq R\norm{\grad\test}{\Lp[\OmegaT]{2}{}} \brac{\norm{\fieldF_k-\fieldF}{\Lp[\OmegaT]{2}{}}+\coeffEL\norm{\fieldEL_k-\fieldEL}{\Lp[\OmegaT]{2}{}}}
     ~~\rightarrow 0.
\end{align*}
Thus, we get for the convection integrals and the electric~drift integrals
\begin{align*}
  \scp{ {\cl[l,k]}\sqbrac{\fieldF_k +\coeffEL\fieldEL_k}-{\cl}\sqbrac{\fieldF+\coeffEL\fieldEL} }{\grad\test}_{\OmegaT} ~~\rightarrow0.
\end{align*}
Finally, for $\Cl_k$ follows from \eqref{eq:existence-conv-2} and Aubin-Lions Lemma the norm-convergence~$\cl[l,k]\rightarrow\cl$ in $\Lp[\OmegaT]{2}{}$.
This leads for the reaction integrals with \Holder\ immediately to
\begin{align*}
 \theta\scp{\Rl(\Cl_k)-\Rl(\Cl)}{\test}_{\OmegaT} \leq \theta\max_lC_{\Rl} \norm{\test}{\Lp[\OmegaT]{2}{}} \sum_l\norm{\cl[l,k]-\cl}{\Lp[\OmegaT]{2}{}}\rightarrow 0.
\end{align*}
In summary, we have shown that $\Cl=\fixOp(\Clold)$ for a arbitrarily chosen subsequence $(\Cl_k)_k$. Therefore, the whole sequence $(\Cl)_k$ converges %(every subsequence converges towards the same limit --> the whole sequence converges)$
and the operator~$\fixOp$ is $\weakstar$-continuous in the sense of equation~\eqref{eq:existence-defWeakstarContinuous}.\\[2.0mm]
%
%====================================================================
% Step 6 - local to global iteration
%====================================================================
\underline{Step~6 -- existence:}~~
A combination of Steps 1 -- 5 shows that we can apply \cref{thm:fpi-theorem}. This yields directly the existence of a solution~$(\fieldEL,\potEL,\fieldF,\press,\Cl)$ on a generally small time interval~$[0,\timeEnd]$.
\par
We now consider for an arbitrary large end time~$\hat{T}$ a time interval $[0,\hat{T}]$, which we decompose with $(K+1)$~time points $0=:T_0 < T_1 < \ldots < T_K :=\hat{T}$ into $K$~subintervals $[T_i,T_{i+1}]$, $i \in \{0,\ldots,K-1\}$.
Furthermore, we suppose that the subintervals $[T_i,T_{i+1}]$ are sufficiently small, such that Steps 1 -- 5 are fulfilled. Thus, a local solution~$(\fieldEL_i,\potEL_i,\fieldF_i,\press_i,\Cl_i)$
exists on $[T_i,T_{i+1}]$ and this solution satisfies the a~priori estimates from \cref{thm:aprioriBounds}. We now carefully check how the constants of the a~priori estimates depend on the
end time~$T_{i+1}$ of the subinterval~$[T_i,T_{i+1}]$. This reveals that the dependency of these constants on $T_{i+1}$ behave as $\exp(T_{i+1})$, which eliminates any possibility of
a blow-up on $[T_i,T_{i+1}]$. Thus, it is admissible to take the partial solution~$\Cl_i$ as the initial value for the $(i+1)$-th solution~$(\fieldEL_{i+1},\potEL_{i+1},\fieldF_{i+1},\press_{i+1},\Cl_{i+1})$.
This leads together with \cref{thm:unique} to a continuation of the solution on the arbitrary large time interval $[0,\hat{T}]$ and consequently to a global solution.
\par
However, we note that this continuation procedure does not lead to solutions on $[0,\infty]$.
\end{proof}
%
%
%%%%%%%%%%%%%%%%%%%%%%%%%%%%%
%%%%%%%%%%%%%%%%%%%%%%%%%%%%%
%  SECTION conclusion       %
%%%%%%%%%%%%%%%%%%%%%%%%%%%%%
%%%%%%%%%%%%%%%%%%%%%%%%%%%%%
\section{Conclusion}
The contribution of this paper was to show the global existence of unique solutions of two-component electrolyte solutions, which are captured by the \dpnp. Here, two-component electrolyte solutions
means that we considered electrolyte solutions, that consist of a neutral solvent and two oppositely charged solutes. In contrast to previous results, we allowed for two oppositely charged solutes with
arbitrary valencies~$\zl[1]>0>\zl[2]$. Most importantly, we successfully established uniform a~priori estimates for the chemical species by using weighted test functions, i.e., instead of the
standard test test functions~$\test=\cl$, we used the weighted test functions~$\test=\abs{\zl}\cl$. By means of this technique we avoided further restrictions such as the electroneutrality constraint of the
volume-additivity constraint. Therefore, the results of this paper apply to general two-component electrolyte solutions, which are captured by the \dpnp. We note, that the a~priori estimates include a uniform
$\Lp[\OmegaT]{\infty}{}$-bound for the charged solutes~$\cl$, which we obtained by the use of Moser's iteration technique. Moreover, the global existence and uniqueness result holds true in two space dimensions
and three space dimensions.
\par
To our best knowledge, in particular for the case of three spatial dimensions, this is the first global existence and uniqueness result for
two-component electrolyte solutions, that firstly are governed by the \dpnp, that secondly include two oppositely charged chemical species with arbitrary valencies, and which thirdly are not subject to further
restrictions such as the electroneutrality constraint, or the volume additivity constraint.
%
%
%%%%%%%%%%%%%%%%%%%%%%%%%%%%%
%%%%%%%%%%%%%%%%%%%%%%%%%%%%%
%  SECTION acknowledgements %
%%%%%%%%%%%%%%%%%%%%%%%%%%%%%
%%%%%%%%%%%%%%%%%%%%%%%%%%%%%
%
\section*{Acknowledgements}
M. Herz was supported by the Elite Network of Bavaria.
%
%
%%%%%%%%%%%%%%%%%%%%%%%%%%%%%
%%%%%%%%%%%%%%%%%%%%%%%%%%%%%
%  SECTION references       %
%%%%%%%%%%%%%%%%%%%%%%%%%%%%%
%%%%%%%%%%%%%%%%%%%%%%%%%%%%%
%
\printbibliography
\end{document}